\documentclass[journal, twocolumn, twoside]{IEEEtran}

\usepackage[all]{xy}
\usepackage{amsmath}
\usepackage{amsthm}
\usepackage{amssymb}
\usepackage{enumerate}
\usepackage{tikz}
\usepackage{pifont} 
\usepackage{mathtools}
\usepackage{algorithm}
\usepackage[noend]{algpseudocode}
\usepackage{mathrsfs}
\usepackage{cases}
\usepackage{stfloats}
\usepackage{blindtext}
\usepackage{url}
\usepackage{tikz-cd}

\newtheorem{theorem}{Theorem}
\newtheorem{lemma}[theorem]{Lemma}

\theoremstyle{definition}
\newtheorem{defi}[theorem]{Definition}
\newtheorem{example}[theorem]{Example}

\newtheorem{cor}[theorem]{Corollary}
\newtheorem{proposition}[theorem]{Proposition}

\newcommand{\KK}{\mathbb{K}} 
\newcommand{\RR}{\mathbb{R}} 
\newcommand{\CC}{\mathbb{C}} 
\newcommand{\PP}{\mathbb{P}} 

\newcommand\norm[1]{\left\lVert#1\right\rVert}
\providecommand{\rank}[1]{\operatorname{rank}(#1)} 
\providecommand{\nnrank}[1]{\operatorname{rank}_{+}(#1)} 

\DeclareMathOperator{\Res}{Res}

\let\O\relax
\DeclareMathOperator{\O}{O}

\allowdisplaybreaks

\begin{document}
\markboth{IEEE Transactions on Information Theory,~Vol.~XX, No.~X, February~2016}{Yang, Comon, and Lim: Uniqueness of Nonnegative Tensor Approximations}

\title{Uniqueness of Nonnegative Tensor Approximations}
 \author{Yang~Qi, Pierre~Comon~\IEEEmembership{Fellow,~IEEE}, and Lek-Heng~Lim\thanks{Manuscript submitted on October 29, 2014; revised November 16, 2015; accepted February 9, 2016. Yang Qi and Pierre Comon are funded by the European Research Council under the European Community's Seventh Framework Program FP7/2007-2013 Grant Agreement no.~320594. Lek-Heng Lim is funded by AFOSR FA9550-13-1-0133, DARPA D15AP00109, NSF IIS 1546413, DMS 1209136, and DMS 1057064.

Yang Qi and Pierre Comon are with CNRS, Gipsa-Lab, University of Grenoble Alpes, F-38000 Grenoble, France (e-mail: yang.qi@gipsa-lab.fr, p.comon@ieee.org). 

Lek-Heng Lim is with the Computational and Applied Mathematics Initiative, Department of Statistics, University of Chicago, 5734 South University Avenue, Chicago, IL 60637, USA (e-mail: lekheng@uchicago.edu).}}
\maketitle



\maketitle

   
\begin{abstract}
We show that for a nonnegative tensor, a best nonnegative rank-$r$ approximation is almost always unique, its best rank-one approximation may always be chosen to be a best nonnegative rank-one approximation, and that the set of nonnegative tensors with non-unique best rank-one approximations form an algebraic hypersurface. We show that the last part holds true more generally for real tensors and thereby determine a polynomial equation so that a real or  nonnegative tensor which does not satisfy this equation is guaranteed to have a unique best rank-one approximation. We also establish an analogue for real or nonnegative symmetric tensors. In addition, we prove a singular vector variant of the Perron--Frobenius Theorem for positive tensors and apply it to show that a best nonnegative rank-$r$ approximation of a positive tensor can never be obtained by deflation. As an aside, we verify that the Euclidean distance (ED) discriminants of the Segre variety and the Veronese variety are hypersurfaces and give defining equations of these ED discriminants.
\end{abstract}
 
\section{Introduction}

Nonnegative tensor decomposition, i.e., a decomposition of a tensor with nonnegative entries (with respect to a fixed choice of bases) into a sum of tensor products of nonnegative vectors, arises in a wide range of applications. These include hyperspectral imaging, spectroscopy,  statistics,  phylogenetics,   data mining, pattern recognition, among other areas; see \cite{LimC09:jchemo,ShasH05:bonn,SmilBG04,ZhanWPP08:josa} and the references therein. One important reason for its prevalence is that such a decomposition  shows  how a joint distribution of discrete random variables decomposes when they are independent  conditional on  a discrete latent random variable \cite{LimC09:jchemo,ZhouBHD14:jasa} --- a ubiquitous model that underlies many applications.
This is in fact one of the simplest Bayesian network
\cite{GSS,Jord04:ss,KollF09}, a local expression of the joint distribution of a set of random variables $x_i$ as
\begin{equation}\label{eq:naive}
p(x_1,\dots,x_d) = \int  \prod_{i=1}^d p(x_i \mid \theta)  \, d\mu_\theta 
\end{equation}
where $\theta$ is some unknown latent random variable. The relation expressed in \eqref{eq:naive} is often called the \textit{naive Bayes hypothesis}. In the case when the random variables $x_1,\dots,x_d$ and the  latent variable $\theta$ take only a finite number of values, the decomposition becomes one of the form
\begin{equation}\label{eq:main}
t_{i_1,\dots, i_d} = \sum\nolimits_{p=1}^r \lambda_r  u_{i_1,p}\cdots u_{i_d,p}.
\end{equation}
One can show \cite{LimC09:jchemo}  that \textit{any} decomposition of a nonnegative tensor of the form in \eqref{eq:main} may, upon normalization by a suitable constant, be regarded as \eqref{eq:naive}, i.e., a marginal decomposition of a joint probability mass function into conditional probabilities under the naive Bayes hypothesis. In the event when the latent variable $\theta$ is not discrete or finite, one may argue that \eqref{eq:main} becomes an approximation with `$ \approx$' in place of `$=$'.

In this article, we investigate several questions regarding nonnegative tensor decompositions and approximations, focusing in particular on uniqueness issues. In Section~\ref{sec:defn}, we define nonnegative tensors in a way that parallels the usual abstract definition of tensors in algebra. We will view them as elements in a tensor product of cones, i.e., tensors in $C_1 \otimes \dots \otimes C_d$ where $C_1,\dots, C_d$ are cones and the tensor product is that of $\mathbb{R}_+$-semimodules  (we write $\mathbb{R}_+ := [0,\infty)$ for the nonnegative reals).  The special case $C_1 = \mathbb{R}_+^{n_1}, \dots, C_d = \mathbb{R}_+^{n_d}$ then reduces to nonnegative tensors.

It has been established in \cite{LimC09:jchemo} that every nonnegative tensor has a best nonnegative rank-$r$ approximation. In Section~\ref{sec:unique} we will show that this best approximation is almost always unique. Furthermore, the set of nonnegative tensors of nonnegative rank  $>r$ that do not have a unique best rank-$r$ approximation form a semialgebraic set contained in a hypersurface. For the special case when $r=1$, we first show in Section~\ref{rank-one-sec} that for a nonnegative tensor, the best nonnegative rank-one and best rank-one approximations coincide. In Section~\ref{sec:gen-unique}, by exploring normalized singular pairs, we find an explicit polynomial expression describing the hypersurface of real (or nonnegative) tensors that admit non-unique best rank-one approximations, which allows one to check whether a given tensor has a unique best rank-one approximation. This polynomial expression also gives a defining equation of the Euclidean distance discriminant of the Segre variety \cite{DrHOSTh15:fcm}. In Section~\ref{sec:symm-unique}, we find results analogous to those in Section~\ref{sec:gen-unique}  for real (or nonnegative) symmetric tensors. We prove an analogue of the Perron--Frobenius theorem for singular values/vectors of positive tensors  in Section~\ref{rank-one-sec} and, among other things, deduce that one cannot obtain a best nonnegative rank-$r$ approximation of a positive tensor by `deflation', i.e., by finding $r$ successive best nonnegative rank-one approximations.

These results would likely shed light on the large number of computational methods for nonnegative matrix factorizations and nonnegative tensor decompositions \cite{AGKM12:pacmstc,Bro98, ChL08:jsc,CZPA09,CoheCC15:spl,FriedH08:oms,Ho08:diss, KP08:jmaa,  KHP14:jgo, KP11:jsc, LeeSeung01:anips,Vava09:jo,ZhouCX12:tsp}.
 
\section{Nonnegative tensors}\label{sec:defn}

A tensor of order $d$ ($d$-tensor for short) may be represented as a $d$-dimensional hypermatrix, i.e., a $d$-dimensional array of (usually) real or complex values. This is a higher-order generalization of the fact that a $2$-tensor, i.e., a linear operator, a bilinear form, or a dyad, can always be represented as a matrix. Such a coordinate representation sometimes hides intrinsic properties --- in particular, this array of coordinates is meaningful only if the bases of underlying vector spaces have been specified in the first place. With this in mind, we prefer to define tensors properly rather than simply regarding them as $d$-dimensional  arrays  of numbers.

The following is the standard definition of tensors. We will see later how we may obtain an analogous definition for nonnegative tensors.
\begin{defi}\label{def:tensor}
Let $V_i$ be a vector space of finite dimension $n_i$ over a field $\KK$, $i=1,\dots, d$, and  let $V_1 \times \dots \times V_d$ be the set of $d$-tuples of vectors. Then the \textit{tensor product} $V=V_1 \otimes \dots \otimes V_d$ is the free linear space spanned by $V_1 \times \dots \times V_d$ quotient 
by the equivalence relation
\begin{multline}\label{eq:rel}
(v_{1},\dots,\alpha v_{i}+\beta v_{i}^{\prime},\dots,v_{d}) \\
\sim \alpha(v_{1},\dots,v_{i},\dots,v_{d})+\beta(v_{1},\dots,v_{i}^{\prime},\dots,v_{d})
\end{multline}
for every $v_i, v_i' \in V_i$, $ \alpha_i, \beta_i \in\KK$, $i=1,\dots,d$. A \textit{tensor} is an element of  $V_1 \otimes \dots \otimes V_d$.
\end{defi}
In particular, \eqref{eq:rel} gives
\begin{equation}\label{eq:scale}
(\alpha_1 v_1,\alpha_2 v_2,\dots,\alpha_d v_d)  = \left(\prod\nolimits_{i=1}^d \alpha_i \right)\, (v_1,v_2,\dots,v_d)
\end{equation}
More details on the definition of tensor spaces may be found in \cite{Como14:spmag,Hack12,Land12,Lim13:HLA}.

A \textit{decomposable tensor} is one of the form $v_1 \otimes \dots \otimes v_d$, $v_i\in V_i$, $i =1,\dots, d$. It represents the equivalence class of tuples up to scaling as in \eqref{eq:scale}, i.e.,
\[
v_1 \otimes \dots \otimes v_d = \left\{ (\alpha_1 v_1,\dots, \alpha_d v_d) : \prod\nolimits_{i=1}^d \alpha_i =1 \right\}.
\]
By \eqref{eq:scale}, it is clear that a decomposable tensor cannot in general be uniquely represented by a $d$-tuple of vectors, what is often called a ``scaling indeterminacy'' in the engineering literature.  When we use the term `unique' in this article, it is implicit that the uniqueness is only up to scaling of this nature. 

From the way a tensor is defined in Definition~\ref{def:tensor}, it is immediate that a nonzero tensor can always be expressed as a  finite sum of nonzero decomposable tensors. When the number of summands is minimal, this decomposition is called a \textit{rank decomposition} (the term ``canonical polyadic'' or \textsc{cp} is often also used) and the number of summands in such a decomposition is called the \textit{rank} of the tensor. In other words, we have the following: 

For every $ T \in V_1 \otimes \dots \otimes V_d$, there exist $ v_{i,p}\in V_i$, $i =1,\dots,d,$ $p = 1,\dots, \rank{T}$, such that 
\begin{equation}\label{CP-def}
T=\sum\nolimits_{p=1}^{\rank{T}} v_{1,p}\otimes \dots \otimes v_{d,p}.
\end{equation}

We present the above material, which is largely standard knowledge, to motivate an analogous construction for real nonnegative tensors. We will first define nonnegative tensors  in a coordinate-dependent manner (i.e.,  depending  on a choice of bases on $V_1,\dots, V_d$), and then in a coordinate-independent manner.

\begin{defi}\label{def:nonneg1}
For each $i=1,\dots, d$, let $V_i$ be a real vector space with $\dim V_i = n_i$. For any fixed choice of basis $\{v_{i,1},\dots, v_{i,n_i}\}$ for $V_i$, we denote by $V_i^+$ the subset of vectors with nonnegative coefficients in $V_i$, i.e.,
\[
V_i^+ =\left\{ \sum\nolimits_{p=1}^{n_i} \alpha_p v_{i,p}  \in V_i : \alpha_1,\dots, \alpha_{n_i} \in \mathbb{R}_+ \right \}.
\]
We will call an element in $V := V_1 \otimes \dots \otimes V_d$ of the form $u_{1} \otimes \dots \otimes u_{d}$ where $u_i \in V_i^+$ for $i =1,\dots,d$, a \textit{nonnegatively decomposable tensor}. The set of \textit{nonnegative tensors} $V^+$ is then the subset of  $V $ defined by
\begin{multline*}
V^+ =\biggl\{  \sum\nolimits_{p=1}^r u_{1,p} \otimes \dots \otimes u_{d,p} \in V : u_{i,p} \in V_i^+,  \\
i = 1,\dots, d,\; p =1,\dots,r, \; r \in \mathbb{N} \biggr\}.
\end{multline*}
\end{defi}

By its definition, every element of $V^+$ has a representation as a finite sum of nonnegatively decomposable tensors.  A decomposition of minimal length then yields the notions of nonnegative tensor rank  and nonnegative tensor rank decomposition. 
\begin{defi}\label{def:nonneg2}
For every $T \in V^+$, there exist $v_{i,p}\in V_i^+$, $i =1,\dots,d$, $p=1,\dots,r$, such that
\begin{equation}\label{eq:ntd}
T=\sum\nolimits_{p=1}^{\nnrank{T}} v_{1,p}\otimes \dots \otimes v_{d,p}
\end{equation}
where
\begin{equation}\label{eq:nr}
\nnrank{T} := \left\{ r : T=\sum\nolimits_{p=1}^{r} v_{1,p}\otimes \dots \otimes v_{d,p} \right\}.
\end{equation}
We will call \eqref{eq:nr} \textit{nonnegative tensor rank} or \textit{nonnegative rank} for short and \eqref{eq:ntd} a \textit{nonnegative rank decomposition} of the nonnegative tensor $T$.
\end{defi}
An obvious property is that $\nnrank{T} \ge \rank{T}$  for any $T \in V^+$.

We now examine an alternative coordinate-free approach to defining nonnegative tensors and nonnegative rank. This approach is also more general, yielding a notion of \textit{conic rank} for a tensor product of any convex cones. We first recall the definition of a tensor product of semimodules. See \cite{Ban13:au} for details on the existence and a construction of such a tensor product.

\begin{defi}\label{defi:tp}
Let $R$ be a commutative semiring and $M, N$ be $R$-semimodules (cf.\ Appendix for the definitions of semirings and semimodules).   A \textit{tensor product} $M \otimes_R N$ of $M$ and $N$ is an $R$-semimodule satisfying the universal property: There is an $R$-bilinear map $\varphi: M \times N \rightarrow M \otimes_R N$ such that given any other $R$-semimodule $S$ together with an $R$-bilinear map $h: M \times N \rightarrow S$, there is a unique $R$-linear map $\tilde{h}: M \otimes_R N \rightarrow S$ satisfying $h=\tilde{h} \circ \varphi$.
\end{defi}

Recall that a \textit{convex cone} $C$ is a subset of a vector space over an ordered field that is closed under linear combinations with nonnegative coefficients, i.e., $\alpha x + \beta y$ belongs to $C$ for all $x, y \in C$ and any nonnegative scalars $\alpha, \beta$.

Since any convex cone $C_i \subset V_i$ is a semimodule over the semiring $\mathbb{R}_+$, we have the unique tensor product of these convex cones  $C_1\otimes\dots\otimes C_d$ as an $\mathbb{R}_+$-semimodule up to isomorphism. More precisely, the \textit{tensor product of cones} $C_{1}\otimes\dots\otimes C_{d}$ is the
quotient monoid $F(C_1,\dots,C_d)/\sim$, where $F(C_1,\dots,C_d)$ is the free monoid generated by all $n$-tuples
$(v_{1},\dots,v_{d})\in C_{1}\times\dots\times C_{d}$, and
$\sim$ is the equivalence relation on $F(C_1,\dots,C_d)$ defined by
\begin{multline*}
(v_{1},\dots,\alpha v_{i}+\beta v_{i}^{\prime},\dots,v_{d}) \\
\sim \alpha(v_{1},\dots,v_{i},\dots,v_{d})+\beta(v_{1},\dots,v_{i}^{\prime},\dots,v_{d})
\end{multline*}
for every $v_{i},v_{i}^{\prime}\in C_{i}$, $\alpha,\beta \in \mathbb{R}_+$, and $i = 1,\dots,d$. The commutative monoid $C_{1}\otimes\dots\otimes C_{d}$ is an $\RR_+$-semimodule. We write $v_{1}\otimes
\dots\otimes v_{d}$ for the equivalence class representing $(v_{1},\dots,v_{d})$ in $F(C_1,\dots,C_d)/ \sim$. 

A \textit{multiconic map}  from $C_{1} \times \dots \times C_{d}$ to a convex cone $C$ is a map $\varphi :C_{1}\times\dots\times C_{d}\to C $ with the property that
\begin{multline*}
\varphi (u_{1},\dots,\alpha v_{i}+\beta w_{i},\dots ,u_{d})\\
=\alpha \varphi(u_{1},\dots,v_{i},\dots ,u_{d})
+\beta \varphi(u_{1},\dots,w_{i},\dots,u_{d})
\end{multline*}
for all $\alpha,\beta\in \mathbb{R}_+$, $i=1,\dots,d$.

The multiconic map $\nu: C_1\times \cdots \times C_m\to C_1\otimes \dots\otimes C_d$ defined by
\[
\nu(v_1,\dots,v_d) = v_1\otimes\dots\otimes v_d \in F(C_1,\dots,C_d)/\sim
\]
and extended nonnegative linearly  to all of $C_{1} \times \dots \times C_{d}$ satisfies the \textit{universal factorization property} often used to define tensor product spaces: If $\varphi $ is a multiconic
map from
$C_1\times \cdots \times C_d$ into a convex cone $C$, then there exists a unique $\RR_+$-linear
map $\psi$ from $C_1\otimes\dots\otimes C_d$ into $C$,
that makes the following diagram
commutative:
\[
\xymatrix{
C_1\times\dots\times C_d\ar[r]^{\nu}\ar[rd]_{\varphi} & C_1\otimes\dots\otimes C_d\ar[d]^{\psi}\\
 & C}
\]
i.e., $\psi\nu =\varphi $. Strictly speaking we should have written  $C_1\otimes_{\mathbb{R}_+} \!\!\dots\otimes_{\mathbb{R}_+}\!\! C_d$ to indicate that the tensor product is one of $\mathbb{R}_+$-semimodules but this is obvious from context. Note that Definition~\ref{defi:tp} is consistent with our earlier definition of nonnegative tensors since $V^+ = V_1^+ \otimes \cdots \otimes  V_d^+ $ as tensor product of cones over $\mathbb{R}_+$.

In \cite{V13:arxiv}, the tensor product of $C_1, \dots, C_d$ is defined to be the convex cone in $V_1 \otimes \cdots \otimes V_d$ formed by $v_1 \otimes \cdots \otimes v_d \in V_1 \otimes \cdots \otimes V_d$, where $v_i \in C_i$, and showed that this tensor product satisfies the above universal factorization property. By the uniqueness of the $\RR_+$-semimodule satisfying the universal property, our construction and the one in \cite{V13:arxiv} are equivalent.

If $C_1 = \mathbb{R}_+^{n_1},\dots, C_d = \mathbb{R}_+^{n_d}$, we may identify
\[
\mathbb{R}_+^{n_1}\otimes \dots \otimes \mathbb{R}_+^{n_d}= \mathbb{R}_+^{n_1 \times \dots \times n_d}
\]
through the interpretation of the tensor product of vectors as a hypermatrix via the Segre outer product
\begin{multline*}  \label{segre:eq}
[v_{1}(1),\dots,v_{1}(n_1)]^{\mathsf{T}}\otimes \cdots \otimes[v_d(1),\dots
,v_d(n_d)]^{\mathsf{T}} \\
= [v_1(i_1) \cdots v_d(i_d)]_{i_1, \dots, i_d=1}^{n_1, \dots, n_d}.
\end{multline*}
Here we write $v(j)$ for the $j$th coordinate of $v \in \mathbb{R}^n$.
 
We note that one may easily extend the notion of nonnegative rank and nonnegative rank decomposition to tensor product of other cones.
\begin{defi}\label{defi:crk}
A tensor $T \in C_1 \otimes \cdots \otimes C_d$ is said to be \textit{decomposable} if $T$ is of the form $u_1 \otimes \cdots \otimes u_d$, where $u_i \in C_i$.
For $T \in C_1 \otimes \cdots \otimes C_d$, the \textit{conic rank} of $T$, denoted by $\nnrank{T}$,  is the minimal value of $r$ such that $T=\sum\nolimits_{p=1}^r u_{1,p} \otimes \dots \otimes u_{d,p}$, where $u_{i,p} \in C_i$, i.e., $T$ is contained in the convex cone generated by $u_{1,1} \otimes \cdots \otimes u_{d,1}, \dots, u_{1,r} \otimes \dots \otimes u_{d,r}$. Such a decomposition will be called a \textit{conic rank decomposition}.
\end{defi}
In the remainder of this paper, we focus our attention on the case $V^+ = V_1^+ \otimes \cdots \otimes  V_d^+$, the convex cone of nonnegative $d$-tensors although we will point out whenever a result holds more generally for arbitrary cones. For any given positive integer $r$, we let
\[
D_r^+ = \{ X\in V_1^+ \otimes \cdots \otimes V_d^+ : \operatorname{rank}_+(X) \leq r \}
\]
denote the set of tensors of nonnegative rank not more than $r$.

\section{Uniqueness of rank decompositions}\label{sec:unique-decomp}

From the standpoints of both identifiability and well-posedness, an important issue is whether a rank decomposition of the form \eqref{CP-def} is unique. It is clear that such decompositions can never be unique when $d=2$, i.e., for matrices. But when $d > 2$, rank decompositions are often unique, which is probably the strongest reason for their utility in applications.
    There are well-known sufficient conditions ensuring uniqueness of rank decomposition \cite{Krus77:laa, SidiB00:jchemo, DoDeLa131:simax, DoDeLa132:simax} and many recent works on the uniqueness of generic tensors of certain ranks \cite{Stra83:laa, ChiaO12:simax, BCO13:am, ChiaOV14:simax}. We highlight three notable results.  
\begin{theorem}[Kruskal]\label{Kruskal-prop}
The rank decomposition of a $d$-tensor $T$ is unique if
\[
\rank{T}\le \frac{1+\sum\nolimits_{i=1}^d (\kappa_i-1)}{2}
\]
where $\kappa_i$ denote the Kruskal rank of the factors $u_{i,1},\dots, u_{i,\rank{T}}$, which is generically equal to the dimension $n_i$ when $n_i \le \rank{T}$.
\end{theorem}

\begin{theorem}[Bocci--Chiantini--Ottaviani]
The rank decomposition of a generic $d$-tensor $T$ of rank-$r$ is unique when
\[
r \leq \frac{\prod_{i=1}^d n_i - (n_1+n_2+n_3-2) \prod_{i=3}^d n_i}{1+\sum_{i=1}^d (n_i-1)}.
\]
\end{theorem}

\begin{theorem}[Chiantini--Ottaviani--Vannieuwenhoven]\label{Chiantini-prop}
The rank decomposition of a generic $d$-tensor $T$ of rank-$r$ is unique when
\[
r < \left\lceil\frac{\prod\nolimits_{i=1}^d n_i}{1+\sum\nolimits_{i=1}^d (n_i-1)}\right\rceil
\]
if $\prod\nolimits_{i=1}^d n_i \leq 15000$, with some exceptional cases.
\end{theorem}
The authors of \cite{ChiaOV14:simax} also strengthened the above result by a prior compression of tensor $T$. The consequence is that the dimensions $n_i$ in Theorem~\ref{Chiantini-prop} may  be replaced by the multilinear rank of $T$,  which allows significant tightening of the upper bound for low multilinear rank tensors.     The maximum $R_{\text{smax}}$ where a generic tensor with rank $\le R_{\text{smax}}$ has a unique rank decomposition has been called the maximum stable rank in \cite{TichPhKol13:tsp}. Theorem~\ref{Chiantini-prop} implies that if $\prod\nolimits_{i=1}^d n_i \leq 15000$, then aside from the exceptional cases, the maximum stable rank is $\left\lceil\prod\nolimits_{i=1}^d n_i/[1+\sum\nolimits_{i=1}^d (n_i-1)]\right\rceil - 1$, which is one less than the (expected) generic rank \cite{Stra83:laa, Lick85:laa, AOP09:tams, CoteBDeLCa09:laa}.  

Nevertheless these results do not apply directly to \textit{nonnegative} decompositions over $\mathbb{R}_+$ (as opposed to decompositions over $\mathbb{C}$) nor to rank-$r$ \textit{approximations} (as opposed to   rank-$r$ decompositions). The purpose of this paper is to provide some of the first results in these directions. In particular, it will be necessary to distinguish between an \textit{exact} nonnegative rank-$r$ decomposition and a best nonnegative rank-$r$ \textit{approximation}. Note that when a best  nonnegative rank-$r$ approximation to a nonnegative tensor $T$ is unique, it means  that 
\begin{equation}\label{eq:appdef}
\min_{\nnrank{X} \le r} \norm{T - X}
\end{equation}
has a unique  minimizer  $X^*$.  The nonnegative rank-$r$ decomposition of  $X^*$  may not however be unique. 

A nonnegative rank decomposition $ X  = \sum\nolimits_{p=1}^r u_{1,p} \otimes \dots \otimes u_{d,p} \in V_1^+ \otimes \cdots \otimes V_d^+$  is said to be \textit{unique} if for any  other  nonnegative rank decomposition $ X  = \sum\nolimits_{p=1}^r v_{1,p} \otimes \dots \otimes v_{d,p}$, there is a permutation $\sigma$ of $\{1, \dots, d\}$ such that $u_{1,p} \otimes \dots \otimes u_{d,p} = v_{1,\sigma(p)} \otimes \dots \otimes v_{d,\sigma(p)}$ for all $p=1, \dots, r$.


\section{Existence and generic uniqueness of rank-$r$ approximations}\label{sec:unique}

Let $V_1,\dots, V_d$ be real vector spaces. Given a nonnegative tensor $T \in V^+$, we consider the best nonnegative rank-$r$ approximations of $T$, where $r$ is less than the nonnegative rank of $T$.
We let
\[
\delta(T) =\inf\nolimits_{X \in D_r^+} \norm{T - X} =\inf\nolimits_{\nnrank{X} \le r} \norm{T - X},
\]
where $\norm{\,\cdot\,}$ is the \textit{Hilbert--Schmidt norm}, i.e., the $l^2$-norm given by the inner product.

Henceforth any unlabelled norm $\norm{\,\cdot\,}$ on $V_1 \otimes \dots \otimes V_d$ will always denote the Hilbert--Schmidt norm. When $d=2$, the Hilbert--Schmidt norm reduces to the Frobenius norm of matrices and when $d=1$, it reduces to the Euclidean norm of vectors. Also, throughout this article, the notation $\langle X,Y \rangle$ will always denote tensor contraction  in all possible indices for $X,Y$ tensors of any order \cite{Lim13:HLA}. When $X$ and $Y$ are of the same order and real, $\langle X, Y\rangle$ reduces to a real inner product and our notation is consistent with the inner product notation; in particular $\langle X, X\rangle = \lVert X\rVert^2$. When $X$ is a $d$-tensor and $Y$ is a $(d-1)$-tensor, $\langle X, Y\rangle$ is a vector --- this is the only other case that will arise in our discussions below. Note however that over $\mathbb{C}$, $\langle \cdot, \cdot \rangle$ is only a symmetric bilinear form and not a complex inner product (which is a sesquilinear form).

\begin{proposition}
Let  $C_i \subseteq V_i^+$ be a closed semialgebraic cone for $i=1,\dots,d$. Then
$D_r^+ = \{ X\in C_1 \otimes \cdots \otimes C_d : \operatorname{rank}_+(X) \leq r \}$ is a closed semialgebraic set.
\end{proposition}
\begin{proof}
It follows from \cite{LimC09:jchemo} that the set is closed and from the Tarski--Seidenberg Theorem \cite{DesiL08:simax} that it is semialgebraic.
\end{proof}

Since $D_r^+$ is a closed set, for any $T \notin D_r^+$, there is some $T^* \in D_r^+$ such that $ \norm{T-T^*} = \delta(T)$. The following result is an analogue of \cite[Theorem~27]{FrieO14:fcm} for nonnegative tensors based on \cite[Corollary~18]{FrieO14:fcm}.
\begin{proposition}\label{genericApproxUniqueness-prop}
Almost every $T \in V^+$ with nonnegative rank $> r$ has a unique best nonnegative rank-$r$ approximation.
\end{proposition}

\begin{proof}
For any $T, T' \in V_1 \otimes \cdots \otimes V_d$, $\lvert \delta(T)-\delta(T')\rvert \leq \norm{T-T'}$, i.e., $\delta$ is Lipschitz and thus differentiable almost everywhere in $V = V_1 \otimes \cdots \otimes V_d$ by Rademacher Theorem.

Consider a general $T \in V^+$. Then in particular $T$ lies in the interior of $V^+$ and there is an open neighborhood of $T$ contained in $V^+$. So $\delta$ is differentiable almost everywhere in $V^+$ as well. Suppose $\delta$ is differentiable at $T \in V^+$. For any $U \in V$, let $\partial \delta^2_T (U)$ be the differential of $\delta^2$ at $T$ along the direction $U$. Since $\norm{T-T^*} = \delta(T)$ we obtain
\begin{align*}
\delta^2(T+tU) &= \delta^2(T) + t \partial \delta^2_T (U) + O(t^2)\\
 &\leq \norm{T+tU-T^*}^2\\
 & = \delta^2(T) + 2t\langle U, T-T^*\rangle + t^2 \norm{U}^2. 
\end{align*}
Therefore, for any $t$, we have $t \partial \delta^2_T(U) \leq 2t\langle U, T-T^*\rangle$, which implies that
\[
\partial \delta^2_T(U) = 2\langle U, T-T^*\rangle.
\]

If $T'$ is another best nonnegative rank-$r$ approximation of $T$, then
\[
2\langle U, T-T^*\rangle = \partial \delta^2_T(U) = 2\langle U, T-T'\rangle,
\]
from which it follows that $\langle T'-T^*, U \rangle = 0$ for any $U$, i.e., $T' = T^*$.
\end{proof}
We note that Proposition~\ref{genericApproxUniqueness-prop} holds more generally for arbitrary closed cones $C_1,\dots, C_d$ in place of $V_1^+,\dots,V_d^+$. 
Our next proposition holds true for arbitrary closed semialgebraic cones $C_1,\dots,C_d$  in place of $V_1^+,\dots,V_d^+$.
\begin{proposition}\label{semi}
The nonnegative tensors satisfying (i) nonnegative rank $> r$, and (ii) do not have a unique best rank-$r$ approximation, form a semialgebraic set that is contained in some hypersurface.
\end{proposition}

\begin{proof}
Observe that $D_r^+$ is the image of the polynomial map
\begin{align*}
\varphi_r \colon (V_1^{+} \times \cdots \times V_d^{+})^r & \to V^+, \\
(u_{1, 1}, \dots, u_{d, 1}, \dots, u_{1, r}, \dots, u_{d, r}) & \mapsto \sum\nolimits_{j=1}^r u_{1, j}\otimes \cdots \otimes u_{d, j}.
\end{align*}
Hence $D_r^+$ is semialgebraic by the Tarski--Seidenberg Theorem \cite{DesiL08:simax} and  the required result follows from \cite[Theorem~3.4]{FrieS16:banach}.
\end{proof}

Now we examine a useful necessary condition for $\sum\nolimits_{p=1}^r T_p$ to be a best rank-$r$ approximation of $T \in V_1\otimes \cdots \otimes V_d$. For a vector $u \in V_i$, we denote by $u(j)$ the $j$th coordinate of $u$, i.e., $u = (u(1), \dots, u(n_i))$, and we will borrow a standard notation from algebraic topology where a hat over a quantity means that quantity is omitted. So for example,
\begin{gather*}
\widehat{u_1} \otimes u_2 \otimes u_3 = u_2 \otimes u_3, \\ u_1 \otimes \widehat{u_2} \otimes u_3 = u_1 \otimes u_3, \\  u_1 \otimes u_2 \otimes \widehat{u_3}  = u_1 \otimes u_2,\\
u_1 \otimes \dots \otimes \widehat{u_i} \otimes \dots \otimes u_d = u_1 \otimes \dots \otimes u_{i-1} \otimes u_{i+1} \otimes \dots \otimes u_d.
\end{gather*}

Let us recall the following well-known fact, which has been used to develop algorithms for nonnegative matrix factorization and nonnegative tensor decomposition.
\begin{lemma}\label{le:rel}
Let $V_1,\dots, V_d$ be real vector spaces and let $T \in V_1 \otimes \dots \otimes V_d$. Let $\rank{T}> r$ and $\lambda  \sum\nolimits_{j=1}^r  T_j$   be a best rank-$r$ approximation, where $T_j = u_{1, j} \otimes \cdots \otimes u_{d, j}$ and $\bigl\|\sum\nolimits_{j=1}^r T_j\bigr\| = 1$. Then for all     $i=1,\dots,d$, and $p=1,\dots, r$,
\begin{multline}\label{wide0-eq}
\langle T, u_{1,p} \otimes \dots \otimes \widehat{u_{i,p}} \otimes \dots \otimes u_{d,p} \rangle \\
= \lambda \left\langle \sum\nolimits_{j=1}^r T_j, u_{1,p} \otimes \dots \otimes \widehat{u_{i,p}} \otimes \dots \otimes u_{d,p} \right\rangle,
\end{multline}
where $\lambda = \langle T, \sum\nolimits_{j=1}^r T_j \rangle$.
  
\end{lemma}

\begin{proof}
Let $L$ denote the line in $V_1 \otimes \cdots \otimes V_d$ spanned by $\sum\nolimits_{j=1}^r v_{1, j} \otimes \cdots \otimes v_{d, j}$, and $L^{\bot}$ denote the orthogonal complement of $L$. Denote the  orthogonal  projection of $T$  onto  $L$ by $\operatorname{Proj}_L (T)$. Then  
\[
\norm{T}^2 = \norm{\operatorname{Proj}_L (T)}^2 + \norm{\operatorname{Proj}_{L^{\bot}}(T)}^2,
\]
and thus
\begin{align*}
\min_{\alpha \ge 0} \biggl\|T &- \alpha \sum\nolimits_{p=1}^r v_{1, p} \otimes \cdots \otimes v_{d, p}\biggr\|^2 \\
&= \norm{T - \operatorname{Proj}_L (T)}^2 = \norm{\operatorname{Proj}_{L^{\bot}}(T)}^2 \\
&= \norm{T}^2 - \norm{\operatorname{Proj}_L (T)}^2.
\end{align*}
So computing
\[
\min_{v_{1,1},\dots,v_{d,r}} \min_{\alpha \ge 0} \norm{T - \alpha \sum\nolimits_{j=1}^r v_{1, j} \otimes \cdots \otimes v_{d, j}}
\]
is equivalent to computing
\[
\max_{v_{1,1},\dots,v_{d,r}}  \operatorname{Proj}_L (T) = \max_{v_{1,1},\dots,v_{d,r}}  \left\langle T, \sum_{j=1}^r v_{1, j} \otimes \cdots \otimes v_{d, j} \right\rangle.
\]

Since $\bigl\|\sum\nolimits_{j=1}^r T_j\bigr\| = 1$, we must have
\[
\left\langle \sum\nolimits_{j=1}^r T_j, u_{1,p} \otimes \dots \otimes \widehat{u_{i,p}} \otimes \dots \otimes u_{d,p} \right\rangle \ne 0
\]
for some $p$. The Jacobian matrix of the hypersurface defined by $\bigl\|\sum\nolimits_{j=1}^r v_{1, j} \otimes \cdots \otimes v_{d, j}\bigr\| = 1$ has constant rank $1$ around $(u_{1, 1}, \dots, u_{d, 1}, \dots, u_{1, r}, \dots, u_{d, r})$, i.e., this real hypersurface is smooth at the point $(u_{1, 1}, \dots, u_{d, 1}, \dots, u_{1, r}, \dots, u_{d, r})$. Hence we may consider the Lagrangian
\begin{multline}\label{wide1-eq}
 \mathcal{L}  = \left\langle T,  \sum\nolimits_{p=1}^r v_{1, p} \otimes \cdots \otimes v_{d, p} \right\rangle\\
 - \lambda \left( \norm{\sum\nolimits_{p=1}^r v_{1, p} \otimes \cdots \otimes v_{d, p}} - 1 \right).
\end{multline}
Setting $\partial\mathcal{L} /\partial v_{i,p} = 0$ at $(u_{1, 1}, \dots, u_{d, 1}, \dots, u_{1, r}, \dots, u_{d, r})$ gives
\begin{multline}\label{wide2-eq}
 \langle T, u_{1,p} \otimes \dots \otimes \widehat{u_{i,p}} \otimes \dots \otimes u_{d,p} \rangle \\
 = \lambda \left\langle \sum\nolimits_{j=1}^r T_j, u_{1,p} \otimes \dots \otimes \widehat{u_{i,p}} \otimes \dots \otimes u_{d,p} \right\rangle
\end{multline}
with $\lambda = \bigl\langle T, \sum\nolimits_{j=1}^r T_j \bigr\rangle$ for all $i=1,\dots,d$, $p=1,\dots, r$.
\end{proof}
Lemma~\ref{le:rel} has a nice geometric interpretation as follows. Let $\widehat{\sigma}_r (\mathbb{P}V_1 \times \cdots \times \mathbb{P}V_d)$ be the cone of the $r$th secant variety of the Segre variety $\mathbb{P}V_1 \times \cdots \times \mathbb{P}V_d$. Suppose $\lambda  \sum\nolimits_{j=1}^r  T_j$ is a smooth point. Then $T - \lambda  \sum\nolimits_{j=1}^r  T_j$ is perpendicular to the tangent space of $\widehat{\sigma}_r (\mathbb{P}V_1 \times \cdots \times \mathbb{P}V_d)$ at $\lambda  \sum\nolimits_{j=1}^r  T_j$.

We presented Lemma~\ref{le:rel} in a concrete affine (as opposed to projective) manner so that there will be no ambiguity when discussing $\lambda$ and  $u_{i,j}$. We will see  later in Definition~\ref{le:ctra} that when $r=1$, these are \textit{normalized singular values} and \textit{normalized singular vector tuples} of $T$.

For a nonnegative tensor $T$ with $\nnrank{T} > r$, we have an inequality in place of the equality in \eqref{wide0-eq}. First we define the \textit{support} of a vector $v \in V$ to be
\[
\operatorname{supp}(v) := \{ i \in \{1,\dots, \dim V \} : v_i \ne 0\}.
\]
\begin{lemma}\label{le:nonnegaperp}
Let $T \in V^+$ with $\nnrank{T} > r$ and $X = \sum_{p=1}^{r'} u_{1,p} \otimes \cdots \otimes u_{d,p}$ be a solution of the optimization problem \eqref{eq:appdef}. Then
\begin{multline}\label{ineq:tan}
\langle T, u_{1,p} \otimes \dots \otimes v_{i,p} \otimes \dots \otimes u_{d,p} \rangle \\
\le \left\langle X, u_{1,p} \otimes \dots \otimes v_{i,p} \otimes \dots \otimes u_{d,p} \right\rangle
\end{multline}
where $v_{i,p} \in V_i^+$, $i = 1, \dots, d$, and $p = 1, \dots, r'$. For each pair $(i, p)$, consider the subspace
\[
\widetilde{V}_{i,p} := \{v \in V_i : \operatorname{supp}(v) \subseteq \operatorname{supp}(u_{i,p}) \}.
\]
Then
\begin{multline}\label{eq:nontan}
\langle T, u_{1,p} \otimes \dots \otimes v_{i,p} \otimes \dots \otimes u_{d,p} \rangle \\
= \left\langle X, u_{1,p} \otimes \dots \otimes v_{i,p} \otimes \dots \otimes u_{d,p} \right\rangle
\end{multline}
for $v_{i,p} \in \widetilde{V}_{i,p}$.
\end{lemma}

\begin{proof}
Fix a pair $(i,p)$ and consider a curve $X(t) = u_{1,p} \otimes \cdots \otimes (u_{i,p}+tv_{i,p}) \otimes \cdots \otimes u_{d,p} + \sum_{j \neq p} u_{1,j} \otimes \cdots \otimes u_{d,j}$, where $v_{i,p} \in V_i^+$. Since for $t \ge 0$, $\norm{T-X(t)}$ achieves a local minimum at $t=0$, i.e., nondecreasing in $[0, \varepsilon)$ for some small $\varepsilon > 0$, the right derivative
\[
\lim_{t \to 0+}\frac{d}{dt} \norm{T-X(t)} \ge 0.
\]
In other words, we have
\begin{multline*}
\langle T, u_{1,p} \otimes \dots \otimes v_{i,p} \otimes \dots \otimes u_{d,p} \rangle \\
\le  \left\langle X, u_{1,p} \otimes \dots \otimes v_{i,p} \otimes \dots \otimes u_{d,p} \right\rangle.
\end{multline*}
In particular, if $v_{i,p} \in \widetilde{V}_{i,p}$, $X(t)$ is nonnegative for $t \in (-\varepsilon, \varepsilon)$, then the local minimality of $\norm{T-X(t)}$ at $0$ implies that
\[
\frac{d}{dt} \norm{T-X(t)} \biggr\rvert_{t=0} = 0,
\]
which gives us
\begin{multline*}
\langle T, u_{1,p} \otimes \dots \otimes v_{i,p} \otimes \dots \otimes u_{d,p} \rangle \\
=\left\langle X, u_{1,p} \otimes \dots \otimes v_{i,p} \otimes \dots \otimes u_{d,p} \right\rangle,
\end{multline*}
as required.
\end{proof}

Recall that  a choice of bases is always implicit when we discuss $V^+$ (cf.\ Definition~\ref{def:nonneg1}) and we may refer to coordinates (or entries) of a nonnegative tensor $T$ without ambiguity.
\begin{lemma}\label{le:nonpos}
Let $T \in V^+$ with $\nnrank{T} > r$ and $X$ be a solution of the optimization problem \eqref{eq:appdef}. Then there exist $i_1, \dots, i_d$ such that the coordinate  $(T-X)_{i_1, \dots, i_d} > 0$.
\end{lemma}

\begin{proof}
Let $X = \sum_{p=1}^{r'} u_{1,p} \otimes \cdots \otimes u_{d,p}$. Suppose $(T-X)_{i_1, \dots, i_d} \le 0$ for all $i_1, \dots, i_d$. Then there is some $p \in \{1, \dots, r'\}$ such that $u_{1,p}(i_1) \cdots u_{d,p}(i_d) > 0$. So
\begin{multline*}
\langle T-X, u_{1,p} \otimes \cdots \otimes u_{d,p} \rangle \\
\le (T-X)_{i_1, \dots, i_d} u_{1,p}(i_1) \cdots u_{d,p}(i_d) < 0,
\end{multline*}
which contradicts \eqref{eq:nontan}. 
\end{proof}

\begin{proposition}\label{prop:nnappdef}
Let $T \in V^+$ with $\nnrank{T} > r$ and $X$ be a solution to the optimization problem \eqref{eq:appdef}. Then $\nnrank{X} = r$.
\end{proposition}

\begin{proof}
Suppose that $\nnrank{X} \le r-1$. By Lemma~\ref{le:nonpos} there is some coordinate $(T-X)_{i_1,\dots,i_d} > 0$. Let $X'$ be the rank-one tensor whose only nonzero coordinate $X'_{i_1,\dots,i_d} = (T-X)_{i_1,\dots,i_d}$. Then $\norm{T-X-X'} < \norm{T-X}$ and $\nnrank{X+X'} \le r$, which contradicts $X$ being a solution of \eqref{eq:appdef}.
\end{proof}

Proposition~\ref{prop:nnappdef} shows that a solution $X$ of \eqref{eq:appdef} indeed has nonnegative rank exactly $r$; so it is in fact appropriate to call $X$ a best nonnegative rank-$r$ approximation of $T$.
 

\section{Rank-one approximations for nonnegative tensors and the Perron--Frobenius theorem}\label{rank-one-sec}
   
We have established in Section~\ref{sec:unique} that a best nonnegative rank-$r$ approximation of a nonnegative tensor is generically unique.  In this section we focus on the case $r =1$ and  find sufficient conditions that guarantee the uniqueness of best nonnegative rank-one approximations. We begin with  the following simple but useful observation: For a nonnegative tensor, a best rank-one approximation can always be chosen to be a best nonnegative rank-one approximation.
 
\begin{theorem}\label{thm:nonneg}
Given $T \in V^+$, let $u_1 \otimes \cdots \otimes u_d \in V_1 \otimes \cdots \otimes V_d$ be a best rank-one approximation of $T$. Then $u_1, \dots, u_d$ can be chosen to be nonnegative, i.e., $u_1 \in V_1^+,\dots,u_d \in V_d^+$.
\end{theorem}

\begin{proof}
Let $T = (T_{i_1, \dots, i_d})$ and $u_i = (u_{i}(1), \dots, u_{i}(n_i))$. Then
\begin{align*}
\lVert T &- u_1 \otimes \cdots \otimes u_d \rVert^2 \\
&= \sum\nolimits_{i_1, \dots, i_d=1}^{n_1,\dots,n_d} \bigl(T_{i_1, \dots, i_d} - u_{1}(i_1)\cdots u_{d}(i_d)\bigr)^2  \\
&\geq  \sum\nolimits_{i_1, \dots, i_d=1}^{n_1,\dots,n_d} \bigl(T_{i_1, \dots, i_d} - \vert u_{1}(i_1) \vert \cdots \vert u_{d}(i_d) \vert \bigr)^2.
\end{align*}
Since $u_1 \otimes \cdots \otimes u_d$ is a best rank-one approximation, we can choose $u_{j}(i_j) = \vert u_{j}(i_j) \vert$, i.e., $u_1 \in V_1^+,\dots,u_d \in V_d^+$. 
\end{proof}

By Theorem~\ref{thm:nonneg}, there is no need to distinguish between a best rank-one and a best nonnegative rank-one approximation of a nonnegative tensor. This allows us to treat best rank-one approximations of a real tensor in a unified way, i.e., we will look for sufficient conditions to ensure a unique best rank-one approximation of a real tensor. Motivated in part by the notion of singular pairs of a tensor \cite{Lim05:camsap} and by the case $r=1$ in Lemma~\ref{le:rel}, we propose the following definition.
\begin{defi}\label{le:ctra}
Let $V_1, \dots, V_d$ be vector spaces over $\KK$ of dimensions $n_1, \dots, n_d$. For $T \in V_1 \otimes \cdots \otimes V_d$, we call $(\lambda, u_1, \dots, u_d) \in \KK \times V_1 \times \dots \times V_d$ a \textit{normalized singular pair} of $T$ if
\begin{equation}\label{def:nsp}
\begin{cases}
\langle T, u_1 \otimes \dots \otimes \widehat{u_i} \otimes \dots \otimes u_d \rangle = \lambda u_i, \\
\langle u_i, u_i \rangle = 1,
\end{cases}
\end{equation}
for all $i =1, \dots, d$. We call $\lambda $ a \textit{normalized singular value} and $(u_1, \dots, u_d)$ is called a \textit{normalized singular vector tuple} corresponding to $\lambda$. If $\KK = \RR$, $\lambda \geq 0$, and $u_i \in V^+_i$, we call $(\lambda, u_1, \dots, u_d)$ a \textit{nonnegative normalized singular pair} of $T$.
\end{defi}

The reader is reminded that  the contraction product $\langle \cdot, \cdot \rangle$ is only an inner product over $\RR$ but not $\CC$. In particular, $\langle u, u \rangle \ne \lVert u \rVert^2 $ over $\CC$. In Definition~\ref{le:ctra} we require that $\langle u_i, u_i \rangle = 1$ instead of $\norm{u_i} = 1$ because $\langle u_i, u_i \rangle = 1$ is an algebraic condition, i.e., it is defined by a polynomial equation. However imposing the condition $\langle u_i, u_i \rangle = 1$ would exclude \textit{isotropic} complex singular vector tuples with $\langle u_i, u_i \rangle = 0$ --- note that over $\CC$ this can happen for $u_i \ne 0$. As such, the following projective variant introduced in \cite{FrieO14:fcm} is useful when we would like to include such isotropic cases.
\begin{defi}\label{def:FOsing}
Let $W_1, \dots, W_d$ be complex vector space. For $T \in W_1\otimes \cdots \otimes W_d$, $([u_1], \dots, [u_d]) \in \mathbb{P}W_1 \times \cdots \times \mathbb{P}W_d$ is called a \textit{projective singular vector tuple} if
\begin{equation}\label{eq:osvt}
\langle T, u_1 \otimes \dots \otimes \widehat{u_i} \otimes \dots \otimes u_d \rangle = \lambda_i u_i
\end{equation}
for some $\lambda_i \in \CC$, $i = 1, \dots, d$.
\end{defi}
The number of projective singular vector tuples of a generic tensor has been calculated in \cite{FrieO14:fcm}. In the sense of \cite{DrHOSTh15:fcm}, this number is the Euclidean distance degree of the Segre variety.

Note that as Definition~\ref{def:FOsing} is over projective spaces, the $\lambda_i$'s are not well-defined complex numbers, and neither is $\prod_{i=1}^d \lambda_i$, but this product corresponds in an appropriate sense to a singular value as we will see next.
 
Definitions~\ref{le:ctra} and \ref{def:FOsing} are related over $\mathbb{C}$ as follows. Suppose $([u_1], \dots, [u_d]) \in \PP W_1\times \cdots \times \PP W_d$ is a projective singular vector tuple. We first choose a representative $(u_1, \dots, u_d)$ of $([u_1], \dots, [u_d])$ that satisfies \eqref{eq:osvt} and has $\norm{u_i} = 1$. Note that we may assume $\prod_{i=1}^d \lambda_i$ to be a nonnegative real number: If $(v_1, \dots, v_d)$ is such that $v_j = e^{i\theta_j} u_j$, then $\langle T, v_1 \otimes \dots \otimes \widehat{v_j} \otimes \dots \otimes v_d \rangle = \mu_j v_j$ and we may choose appropriate $\theta_1, \dots, \theta_d$ so that
\[
\prod\nolimits_{i=1}^d \mu_i = e^{i(d-2)(\theta_1 + \dots+ \theta_d)} \prod\nolimits_{i=1}^d \lambda_i \in \RR_+.
\]
For a  nonnegative $\prod_{i=1}^d \lambda_i$,
\[
\lambda := \left( \prod\nolimits_{i=1}^d \lambda_i \right)^{1/d}
\]
is `almost' a normalized singular value of $T$ with corresponding normalized singular vector tuple $(u_1, \dots, u_d)$ --- `almost' because the condition $\langle u_i, u_i \rangle = 1$  in Definition~\ref{le:ctra}  has to be replaced by $\norm{u_i} = 1$.



It has been shown in \cite{FrieO14:fcm} that a generic $T$ does not have a zero singular value nor  a projective singular vector tuple $([u_1], \dots, [u_d])$ such that $\langle u_i, u_i \rangle = 0$ for some $i$. Thus, for a generic $T$, both definitions above are equivalent. We may use the two definitions interchangeably depending on the situation. In this article, we will mainly consider the normalized singular pairs of a tensor as defined in Definition~\ref{le:ctra}.

The next three results give an analogue of the tensorial Perron--Frobenius Theorem \cite{ ChanPZ08:cms, FrieGH13:laa, Lim05:camsap, YangYang10:jmaa} for nonnegative normalized singular pairs (as opposed to nonnegative eigenpairs \cite{Lim05:camsap}). The proof of Lemma~\ref{exi} in particular will require the $l^1$-norm. Again recall that  a choice of bases is always implicit when we discuss $V^+$ (cf.\ Definition~\ref{def:nonneg1}) and the $l^1$-norm is with respect to this choice of bases.
\begin{lemma}[Existence]\label{exi}
A nonnegative tensor $T \in V^+$ has at least one nonnegative normalized singular pair.
\end{lemma}

\begin{proof}
Consider the compact convex set
\[
D = \left\{ (u_1, \dots, u_d) \in V_1^+ \times \dots \times V_d^+ : \sum\nolimits_{i=1}^d \| u_{i} \|_1 = 1 \right\}.
\]
If $ \sum\nolimits_{i=1}^d \|\langle T, u_1 \otimes \dots \otimes \widehat{u_i} \otimes \dots \otimes u_d \rangle \|_1 = 0$ for some $(u_1, \dots, u_d)$, then $\langle T, u_1 \otimes \dots \otimes \widehat{u_i} \otimes \dots \otimes u_d \rangle = 0$ for all $i$, which implies that $\lambda = 0$. On the other hand,  if  $ \sum\nolimits_{i=1}^d \|\langle T, u_1 \otimes \dots \otimes \widehat{u_i} \otimes \dots \otimes u_d \rangle \|_1 > 0$, we define the map $\psi \colon D \to D$ by\par
\vspace*{-2ex}
{\footnotesize
\begin{align*}
&\psi (u_1, \dots, u_d) \\
&= \biggl(  \frac{\langle T, u_2 \otimes \cdots \otimes u_d \rangle}{ \sum\nolimits_{i=1}^d \| \langle T, u_1 \otimes \dots \otimes \widehat{u_i} \otimes \dots \otimes u_d \rangle\|_1},\dots \\
&\qquad \qquad \qquad \qquad \dots, \frac{\langle T, u_1 \otimes \cdots \otimes u_{d-1} \rangle}{ \sum\nolimits_{i=1}^d \| \langle T, u_1 \otimes \dots \otimes \widehat{u_i} \otimes \dots \otimes u_d \rangle \|_1} \biggr).
\end{align*}}%
Note that each term $\langle T, u_1 \otimes \dots \otimes \widehat{u_i} \otimes \dots \otimes u_d \rangle$ in the denominator is the contraction of a $d$-tensor with a $(d-1)$-tensor and therefore the result is a vector. We then normalize by the sum of the $l^1$-norms of these vectors so that $\|\psi\|_1=1$.  

By Brouwer's Fixed Point Theorem, there is some $u_1 \otimes \cdots \otimes u_d$ such that $\langle T, u_1 \otimes \dots \otimes \widehat{u_i} \otimes \dots \otimes u_d \rangle = \lambda u_i$ where
\[
\lambda = \sum\nolimits_{i=1}^d \|\langle T, u_1 \otimes \dots \otimes \widehat{u_i} \otimes \dots \otimes u_d \rangle \|_1.
\]
Since $\langle T, u_1 \otimes \cdots \otimes u_d \rangle = \lambda  \norm{u_i}^2$ for $i = 1, \dots, d$, $\norm{u_1} = \cdots = \norm{u_d}$. Let $u'_i = u_i /\norm{u_i}$ and $\lambda' = \langle T, u'_1 \otimes \cdots \otimes u'_d \rangle$. Then $(\lambda', u'_1, \dots, u'_d)$ is a nonnegative normalized singular pair.
\end{proof}
One of our reviewers has pointed out to us that Lemma~\ref{exi} may also be obtained from Lemma~\ref{le:rel} and Theorem~\ref{thm:nonneg}.

\begin{defi}\label{def:pos}
We say that a tensor $T \in V^+$ is \textit{positive} if all its coordinates (with respect to the implicit choice of bases when we specify $V^+$, cf.\ Definition~\ref{def:nonneg1})  are positive.
\end{defi}
\begin{lemma}[Positivity]\label{pos}
If $T$ is positive, then $T$ has a positive normalized singular pair $(\lambda, u_1, \dots, u_d)$ with $\lambda > 0$.
\end{lemma}

\begin{proof}
By Lemma~\ref{exi}, $T$ has a nonnegative normalized singular pair $(\lambda, u_1, \dots, u_d)$. Suppose that a choice of bases has been fixed for $V_1,\dots,V_d$. We let $v_i(j)$ denote the $j$th coordinate of a vector $v_i \in V_i$, $j = 1,\dots,n_i$.
Let
\[
\alpha = \min\{ u_{i}(j) : i =1 ,\dots, d, \; j\in \operatorname{supp}(u_i)\}.
\]
For any $i$ and $j$,
\begin{align*}
\lambda u_{i}(j) &= \langle T, u_1 \otimes \dots \otimes \widehat{u_i} \otimes \dots \otimes u_d \rangle(j) \\
& \geq \alpha^{d-1} \sum\nolimits_{k_j \in \operatorname{supp}(u_j)} T_{k_1\dots k_{i-1} j k_{i+1}\dots k_d} > 0,
\end{align*}
implying that $\lambda$ and all coordinates of $u_i$ are positive.
\end{proof}

We recall the definition of \textit{spectral norm} for a tensor, which is known  \cite{HL13:jacm}  to be  NP-hard to compute or even approximate.
\begin{defi}
For $T \in V_1 \otimes \cdots \otimes V_d$ over $\RR$, let $\lVert T \rVert_\sigma := \max \{ \vert \langle T, u_1 \otimes \cdots \otimes u_d \rangle \vert: \norm{u_1} = \dots = \norm{u_d} = 1 \}$ be the \textit{spectral norm} of $T$.
\end{defi}
 
We may deduce the following from \cite[Theorem~20]{FrieO14:fcm} and Lemma~\ref{le:rel}.
\begin{cor}[Generic Uniqueness]\label{uni}
A general real tensor $T$ has a unique normalized singular pair $(\lambda, u_1, \dots, u_d)$ with $\lambda = \lVert T \rVert_\sigma$.
\end{cor}

 

The relation between best rank-$r$ and best rank-one approximations of a matrix over $\mathbb{R}$ or $\mathbb{C}$ is well-known: A best rank-$r$ approximation can be obtained from $r$  successive best rank-one approximations --- a consequence of the Eckart--Young Theorem.  It has been shown in  \cite{StegC10:laa} that this `deflation procedure' does not work for real or complex $d$-tensors of order $d > 2$. In fact, more recently, it has been shown in \cite{VanNVM14:jmaa} that the property almost never holds when $d > 2$.

We will see here that the `deflatability' property does not  hold for nonnegative tensor rank either.
   
\begin{proposition}\label{sey}
A best nonnegative rank-$r$ approximation of a positive tensor with nonnegative rank $> r$ cannot be obtained by a sequence of best nonnegative rank-one approximations.
\end{proposition}

\begin{proof}
It suffices to show that a best nonnegative rank-$2$ approximation cannot be obtained by two best nonnegative rank-one approximations. Let $T \in V^+$ be a positive tensor with $\nnrank{T} > 2$. Suppose $u_1 \otimes \cdots \otimes u_d$ is a best rank-one approximation of $T$, and $u_1 \otimes \cdots \otimes u_d + v_1 \otimes \cdots \otimes v_d$ is a best nonnegative rank-$2$ approximation of $T$. By the proof of Lemma~\ref{pos}, $u_k > 0$ for all $k=1,\dots, d$, then by Lemma~\ref{le:nonnegaperp}, we have
\begin{gather*}
\left\langle T - u_1 \otimes \cdots \otimes u_d, \, u_1 \otimes \cdots \otimes u_d  \right\rangle =  0,\\
\left\langle T - u_1 \otimes \cdots \otimes u_d - v_1 \otimes \cdots \otimes v_d, \,  u_1 \otimes \cdots \otimes u_d  \right\rangle =  0.
\end{gather*}
We subtract the second equation from the first to get
\[
\left\langle v_1 \otimes \cdots \otimes v_d, \,  u_1 \otimes \cdots \otimes u_d  \right\rangle = 0,
\]
which contradicts the non-negativity of each $v_k$ and the positivity of each $u_k$.
\end{proof}
 
Following \cite{VanNVM14:jmaa}, we say that a tensor $T \in V^+$ with nonnegative rank $s$ admits a \textit{Schmidt--Eckart--Young decomposition} if it can be written as a linear combination of nonnegatively decomposable tensors $T = \sum\nolimits_{p=1}^s u_{1,p} \otimes \dots \otimes u_{d,p}$, and such that $\sum\nolimits_{p=1}^r u_{1,p} \otimes \dots \otimes u_{d,p}$ is a best nonnegative rank-$r$ approximation of $T$ for all $ r=1, \dots, s$. Proposition~\ref{sey} shows that a general nonnegative tensor does not admit a Schmidt--Eckart--Young decomposition.

We point out  that methods in \cite{dSCdA15:arxiv, PhTC151:tsp, PhTC152:tsp} (for real/complex) \cite{CZPA09, ZhouCX12:tsp, KHP14:jgo} (nonnegative) rely on deflation.

\section{Uniqueness of best rank-one approximations for real symmetric tensors}\label{sec:symm-unique}

Not every tensor has a unique best rank-one approximation \cite[Proposition~1]{StegC10:laa}. For example, the symmetric $3$-tensor $x \otimes x \otimes x + y \otimes y \otimes y$, where $x$ and $y$ are orthonormal, has two best rank-one approximations: $x\otimes x \otimes x$ and $y \otimes y \otimes y$. It is known that a best rank-one approximation of a symmetric tensor can be chosen to be symmetric over $\mathbb{R}$ and $\mathbb{C}$ \cite{Banach,FrieS13:fmc}. In this section we study various properties of the set of symmetric tensors that do not have unique  best symmetric rank-one approximations. Before we get to these we will have to first introduce analogues/generalizations of  \textit{eigenpairs} and \textit{characteristic polynomials}  for higher-order symmetric tensors.

In the following, for a real or complex vector space $V$, $\mathsf{S}^d(V)$ denotes the symmetric $d$-tensors over $V$. For any $u \in V$, we write $u^{\otimes d} = u \otimes \dots \otimes u \in  \mathsf{S}^d(V)$ for the $d$-fold tensor product of $u$ with itself.

Let $V^*$ be the dual space of $V$. For any group $G$ acting on $V$, $G$ also acts naturally on $\mathsf{S}^d(V)$ and $\mathsf{S}^d(V^*)$ such that
\[
\langle S, T \rangle = \langle g \cdot S, g \cdot T \rangle
\]
for all $g \in G$, $T \in \mathsf{S}^d(V)$, and $S \in \mathsf{S}^d(V^*)$. If we fix an inner product $(\cdot, \cdot)$ on $V$, then $V$ becomes self dual and we may identify $V^*=V$. In which case $\langle \cdot, \cdot \rangle$ may be regarded  the inner product on $\mathsf{S}^d(V)$ defined by
\[
\langle u^{\otimes d}, v^{\otimes d} \rangle := (u, v)^d 
\]
and extended linearly to any $S, T \in  \mathsf{S}^d(V)$ (since any element of $ \mathsf{S}^d(V)$ may be expressed as a linear combination of $u^{\otimes d}$'s \cite{CGLM}). The inner product $\langle \cdot, \cdot \rangle$ is clearly invariant under the group that preserves the inner product  $(\cdot, \cdot)$. In particular, if $V = \mathbb{R}^n$, then $\langle \cdot, \cdot \rangle$  is invariant under the orthogonal group\footnote{Henceforth we assume that our vector spaces are equipped with inner products and we write $\O(n)$ for the group that preserves the inner product.} $\O(n)$.

The following definition of symmetric tensor eigenpairs is based on \cite{CS13:laa, Lim05:camsap, Qi05:jsymbcomp}.
\begin{defi}\label{def:eigpair}
For $T \in \mathsf{S}^d(V)$ over $\CC$, $(\lambda, u) \in \CC \times V$ is called a \textit{normalized eigenpair} of $T$ if
\[
\begin{cases}
\langle T, u^{\otimes (d-1)}\rangle = \lambda u,\\
\langle u, u \rangle = 1.
\end{cases}
\]
$\lambda$ is the \textit{normalized eigenvalue} and $v$ the corresponding \textit{normalized eigenvector} of $T$.
Two normalized eigenpairs $(\lambda, u)$ and $(\mu, v)$ of $T$ are \textit{equivalent} if $(\lambda, u) = (\mu, v)$ or if $(-1)^{d-2}\lambda = \mu$ and $u = -v$. A normalized eigenvalue $\lambda$ is said to be \textit{simple} if it has only one corresponding normalized eigenvector up to equivalence.
\end{defi}
The number of eigenpairs of a tensor over $\CC$ has been determined in \cite{CS13:laa, OO13:jsymcomp}; one may view this as the ED degree of the Veronese variety \cite{DrHOSTh15:fcm}.
Definition~\ref{def:eigpair} also applies to a real vector space $V$. In this case, normalized eigenpairs of  $T \in \mathsf{S}^d(V)$ are invariant $\O(n)$.

It is easy to see that for a symmetric tensor $T \in   \mathsf{S}^d(V)$, the spectral norm $\lVert T \rVert_\sigma$ is the largest eigenvalue of $T$ in absolute value. 
Let $\mathbb{S}^{n-1} $ denote the unit sphere in $\RR^n$.   The subset $\{ u \in \mathbb{S}^{n-1}: \langle T, u^{\otimes d} \rangle = \lVert T \rVert_\sigma \}$ is non-empty and closed in $\mathbb{S}^{n-1}$ and invariant under $\O(n)$.

To introduce the characteristic polynomial of a symmetric tensor, we first recall the definition and some basic properties of the \textit{multipolynomial resultant} \cite{GKZ94, CLO05}. For any given $n+1$ homogeneous polynomials $F_0, \dots, F_n \in \CC[x_0, \dots, x_n]$ with positive total degrees $d_0, \dots, d_n$, let $F_i = \sum_{\vert \alpha \vert = d_i} c_{i, \alpha} x_0^{\alpha_0} \cdots x_n^{\alpha_n}$, where $\alpha = (\alpha_0,\dots,\alpha_n)$ and $\lvert \alpha \rvert = \alpha_0 + \cdots + \alpha_n$. We will associate each pair $(i, \alpha)$ with a variable $u_{i, \alpha}$. Now given a polynomial $P$ in the variables $u_{i, \alpha}$ where $i=0,\dots,n$, and $\lvert \alpha \rvert \in \{d_0,\dots,d_n\}$, we denote by $P(F_0, \dots, F_n)$ the result obtained by substituting each $u_{i, \alpha}$ in $P$ with $c_{i, \alpha}$. The following is a classical result in invariant theory \cite{GKZ94, CLO05}.
\begin{theorem}
There is a unique polynomial $\Res$ with integer coefficients in the variables $u_{i, \alpha}$ where $i=0,\dots,n$, and $\lvert \alpha \rvert \in \{d_0,\dots,d_n\}$,  that has the following properties:
\begin{enumerate}[\upshape (i)]
\item $F_0 = \cdots = F_n = 0$ have a nonzero solution over $\CC$ if and only if $\Res{(F_0, \dots, F_n)} = 0$.
\item $\Res{(x_0^{d_0}, \dots, x_n^{d_n})} = 1$.
\item $\Res$ is irreducible over $\CC$.
\end{enumerate}
\end{theorem}

\begin{defi}
$\Res{(F_0, \dots, F_n)} \in \CC$ is called the \textit{resultant} of the polynomials $F_0, \dots, F_n$. Often we will also say that it is the resultant of the system of polynomial equations $F_0 =0, \dots, F_n = 0$.
\end{defi}

The following definition was first proposed in \cite{Qi07:jmma} and called an $E$-characteristic polynomial.
\begin{defi}\label{def:Echar}
The \textit{characteristic polynomial} of a symmetric tensor $T$ is the resultant $\psi_T(\lambda)$ of  the following systems of polynomial equations in $n+1$ variables $u$ and $x$ (note that $u$ has $n$ entries).
\begin{enumerate}[\upshape (i)]
\item 
For $T \in \mathsf{S}^{2d-1}(V)$,
\[
\langle T, u^{\otimes (d-1)} \rangle - \lambda x^{d-2} u = 0 \quad \text{and}\quad x^2 - \langle u, u \rangle = 0.
\]

\item For $T \in \mathsf{S}^{2d}(V)$,
\[
\langle T, u^{\otimes (2d-1)} \rangle - \lambda \langle u, u \rangle^{d-1}u  = 0.
\]
\end{enumerate}
\end{defi}
Note that we regard $\lambda$ as a parameter and not one of the variables. One may show that the resultant $\psi_T(\lambda)$ is a (univariate) polynomial in $\lambda$.

In the following, for $u,v,w \in V$, we write
\begin{multline*}
u \odot v \odot w := \frac{1}{6}(u \otimes v \otimes w + u \otimes w \otimes v + v \otimes u \otimes w\\
 + v \otimes w \otimes u + w \otimes v \otimes u + w \otimes u \otimes v)
\end{multline*}
for the \textit{symmetric tensor product} \cite{CGLM}. Note that $u \odot v \odot w = v \odot u \odot w = \dots = w \odot v \odot u$, i.e., symmetric tensor product is independent of order and in particular $u \odot v \odot w  \in \mathsf{S}^3(V)$. It is easy to extend this to arbitrary order
\[
u_1 \odot \cdots \odot u_d = \frac{1}{d!} \sum\nolimits_{\tau \in \mathfrak{S}_d} u_{\tau(1)} \otimes \cdots \otimes u_{\tau(d)} \in \mathsf{S}^d(V).
\]
For $u \in V$, we may write $u^{\odot d} = u \odot \cdots \odot u$ for the $d$-fold symmetric tensor product of $u$ with itself but  we clearly always have
\[
u^{\otimes d} = u^{\odot d}.
\]

\begin{proposition}\label{thm:3symhy}
Let $V$ be a real vector space of dimension $n$. Let  $\rho = \lVert T \rVert_\sigma$ and define
\[
H_\rho := \{ T \in \mathsf{S}^d(V) :  \rho \; \text{is not simple}\}. 
\]
Then $H_{\rho}$ is an algebraic hypersurface in $\mathsf{S}^d(V)$. 
\end{proposition}

\begin{proof}
For notation convenience, we prove the result for $d = 3$; extending to $d > 3$ is straightforward. Let $T \in \mathsf{S}^3(V)$. Suppose $T \in H_{\rho}$, i.e., there exist $u \neq v\in V$ with $\norm{u} = \norm{v} = 1$ such that
\[
\langle T, u^{\otimes 2}\rangle = \rho u, \qquad \langle T, v^{\otimes 2}\rangle = \rho v.
\]

Let $u_1 := u$ and extend it to $\{u_1, \dots, u_n \}$, an orthonormal basis of $V$. By an action of the orthogonal group $\O(n)$ on $V$, we may assume that $v = u_1\cos \theta  + u_2 \sin \theta $ for some $\theta \in (0,\pi)$. Let $T_{ijk} := \langle T, u_i \odot u_j \odot u_k\rangle $. Then
\begin{numcases}{}
T_{111} = \rho,  \nonumber \\
T_{i11} = 0, \\ \label{eq:sympara3}
T_{111}\cos^2 \theta + T_{122}\sin^2 \theta = T_{111} \cos \theta, \\ \label{eq:sympara4}
2T_{122} \sin \theta \cos \theta + T_{222} \sin^2 \theta = T_{111} \sin \theta,  \nonumber \\ 
2T_{j12} \cos \theta + T_{j22} \sin \theta = 0,  \nonumber 
\end{numcases}
for $i \neq 1$ and $j > 2$.

By eliminating $\theta$, we may obtain equations in $T_{ijk}$'s. For example, \eqref{eq:sympara3} implies $\cos \theta = 1$ or $(T_{111} - T_{122}) \cos \theta = T_{122}$, and \eqref{eq:sympara4} implies $\sin \theta = 0$ or $2T_{122}\cos \theta + T_{222} \sin \theta = T_{111}$. Since $\sin^2 \theta + \cos^2 \theta = 1$ and $\theta \neq 0$ or  $\pi$, we have\par
\vspace*{-2ex}
{\footnotesize
\begin{equation}\label{eq:symuni}
\begin{cases}
[T_{111}(T_{111} - T_{122}) - 2T^2_{122}]^2 + T^2_{222}T^2_{122} = T^2_{222}(T_{111} - T_{122})^2, \\
(T_{111}T_{122} + 2T^2_{122} - T^2_{111})T_{j22} = 2T_{j12}T_{222}(T_{111} - T_{122}).
\end{cases}
\end{equation}}%

Let $J := \{ (T, [u_1, \dots, u_n]) \in \mathsf{S}^3(V) \times \O(n) : T_{ijk}$ satisfies \eqref{eq:symuni}$\}$. Consider the projections
\begin{equation}
\xymatrix{
& J \ar[ld]^{\pi_{1}} \ar[rd]_{\pi_{2}} &\\
\mathsf{S}^3(V) & & \O(n) }
\end{equation}
where $\pi_1 (T, [u_1, \dots, u_n]) = T$ and $\pi_2(T, [u_1, \dots, u_n] ) = [u_1, \dots, u_n]$. By \cite{Qi07:jmma}, $\rho$ is a root of the $E$-characteristic polynomial $\psi_T(\lambda)$ of $T$. So $\rho$ and any of its corresponding normalized eigenvectors must depend algebraically on $T$, implying that $J$ is a variety in $\mathsf{S}^3(V) \times \O(n)$. Hence $T$ has more than one inequivalent normalized eigenvectors corresponding to $\rho$ if and only if $T$ is in the image of $\pi_1$, i.e., $H_{\rho} = \pi_1 (J)$.

Now define $T' \in \mathsf{S}^3(V)$ by
\[
T'_{111} = 1, \quad  T'_{122} = 2\sqrt{3}-3, \quad T'_{222} = 6\sqrt{3}-10,
\]
and set all other terms $T'_{ijk} = 0$. Then $T'$ has two normalized eigenvectors corresponding to its normalized eigenvalue $\rho = \lVert T' \rVert_\sigma = 1$. Hence $T' \in \pi_1(J)$. Since $T'$ has a finite number of eigenvectors, a generic $T \in \pi_1(J)$ must also have a finite number of eigenvectors by semicontinuity. Hence $\dim \pi_1^{-1}(T) = \dim \O(n-2)$ for a generic $T \in \pi_1(J)$. So $\dim H_{\rho} = \dim \pi_1(J) = \dim J - \dim \O(n-2) = \dim J -  (n-2)(n-3)/2$.

Since $\pi_2$ is a dominant morphism, and the dimension of a generic fiber $\pi_2^{-1}([u_1, \dots, u_n])$ is $\dim \mathsf{S}^3(V) - 2(n-1)$, we deduce that $\dim J = \dim \mathsf{S}^3(V) - 2(n-1) + \dim \O(n)$. Therefore $\dim H_{\rho} = \dim \mathsf{S}^3(V) - 1$, i.e., $H_{\rho}$ is a hypersurface.
\end{proof}
\begin{figure*}[b!]
\begin{equation}\label{eq:G}
G = \left[ \begin{smallmatrix}
T_{111} & T_{122} & 0 & 2T_{112} & -\lambda & 0 \\
T_{112} & T_{222} & 0 & 2T_{122} & 0 & -\lambda \\
1 & 1 & -1 & 0 & 0 & 0 \\
12T_{122} & 4T_{111}\lambda - 8T_{122}\lambda & 4T_{111}\lambda + 4T_{122}\lambda & 8T_{222}\lambda - 16T_{112}\lambda & 16T_{112}^2 - 4\lambda^2 - 16T_{111}T_{122} & 8T_{112} T_{122} - 8T_{111} T_{222} \\
4T_{222}\lambda - 8T_{112}\lambda & 12T_{112}\lambda & 4T_{112}\lambda + 4T_{222}\lambda & 8 T_{111}\lambda - 16T_{122}\lambda & 8T_{112}T_{122} - 8T_{111}T_{222} & 16T_{122}^2 - 4\lambda^2 - 16T_{112}T_{222} \\
8T_{112}^2 - 8T_{111}T_{122} - 2\lambda^2 & 8T_{122}^2- 8T_{222} T_{112} - 2\lambda^2 & - 6\lambda^2 & 8T_{112}T_{122} - 8T_{111}T_{222} & 8T_{122}\lambda + 8T_{111}\lambda & 8T_{112}\lambda + 8T_{222}\lambda \end{smallmatrix} \right]
\end{equation}
\end{figure*}

Let $V$ be a real vector space. We specify a choice of basis on $V$ and define the set of nonnegative symmetric tensors to be
\[
\mathsf{S}^d(V^+) := \mathsf{S}^d(V) \cap ( V^{\otimes d})^+.
\]
Recall also Definition~\ref{def:pos}.
\begin{cor}
Let $T \in \mathsf{S}^3(V^+)$ be positive. Let $u \in V$ be such that
$\langle T, u^{\otimes 3}\rangle = \rho = \lVert T \rVert_\sigma$
and
\[
\sigma_2 := \min \{ \vert \langle T, u \odot v \odot v\rangle \vert: \langle u, v\rangle =0, \; \norm{v}=1 \}.
\]
If $\sigma_2 \geq \rho/2 $, then $T$ has a unique best nonnegative symmetric  rank-one approximation.
\end{cor}

\begin{proof}
By Lemma~\ref{le:rel}, suppose there exist $v \neq u$ such that $\norm{v} = 1$, $\langle u, v\rangle =0$, and $\langle T, (u\cos \theta + v\sin \theta)^{\otimes 3}\rangle = \rho$ for some $0 < \theta \leq \pi$. Then by Lemma~\ref{pos}, we must in fact have $0 < \theta < \pi/2$. By \eqref{eq:sympara3}, $\langle T, u \odot v \odot v\rangle = \frac{\cos \theta}{1 + \cos \theta} \rho$. Since $0 < \frac{\cos \theta}{1 + \cos \theta} < \frac{1}{2}$ when $0 < \theta < \pi/2$, we get $0 < \langle T, u \odot v \odot v\rangle < \rho/2$, which contradicts $\sigma_2 \ge \rho/2$.
\end{proof}

Let $V$ be a real vector space of dimension $n$ and $W = V \otimes_{\RR} \CC$ be its complexification. A generic $T \in \mathsf{S}^d(W)$ has distinct eigenvalues \cite{CS13:laa}, so the resultant of the polynomial $\psi_T$ and its derivative $\psi'_T$, denoted by $D_{\operatorname{eig}}(T)$, is a nonzero polynomial on $\mathsf{S}^d(W)$ called the \textit{eigen discriminant}. The equation $D_{\operatorname{eig}}(T) = 0$ defines the complex hypersurface $H_{\operatorname{disc}}$ consisting of tensors $T \in \mathsf{S}^d(W)$ that do not have simple normalized eigenpairs. For $T \in \mathsf{S}^d(V)$, the hypersurface $H_{\rho}$ in Proposition~\ref{thm:3symhy} is a union of some components of the real points of $H_{\operatorname{disc}}$. In fact, if we replace $\rho = \lVert T \rVert_\sigma$ by any real normalized eigenvalue $\mu$ of $T$ in the proof of Proposition~\ref{thm:3symhy}, we may show that the subset of symmetric tensors whose normalized eigenvalues are not all simple is a finite union of real algebraic hypersurfaces, and these hypersurfaces are the real points of $H_{\operatorname{disc}}$. We summarize this discussion as follows.

\begin{theorem}\label{EDsym}
$D_{\operatorname{eig}}(T) = 0$ is a defining equation of the hypersurface
\[
H_{\operatorname{disc}} := \{ T \in \mathsf{S}^d(W) :  T \; \text{has a non-simple eigenvalue}\}. 
\]
\end{theorem}

For $T \in \mathsf{S}^d(V)$, if  $D_{\operatorname{eig}}(T) \neq 0$, then by definition, either (i) there is a unique eigenvector $v_{\lambda}$ corresponding to each eigenvalue $\lambda$ of $T$ when $d$ is odd, or (ii) there are two eigenvectors $\pm v_{\lambda}$ corresponding to each eigenvalue $\lambda$ of $T$ when $d$ is even. Hence we have the following.
\begin{cor}\label{cor:EDsym}
Let $T \in \mathsf{S}^d(V)$. If $D_{\operatorname{eig}}(T) \neq 0$, then $T$ has a unique best  symmetric rank-one approximation.
\end{cor}
We deduce the following analogue for nonnegative tensors from Banach's Theorem that the best rank-one approximation of a symmetric tensor can be chosen to be symmetric \cite{Banach,FrieS13:fmc},  Theorem~\ref{thm:nonneg}, and Corollary~\ref{cor:EDsym}.
\begin{cor}\label{cor:EDnonneg}
Let $T \in \mathsf{S}^d(V^+)$. If $D_{\operatorname{eig}}(T) \neq 0$, then $T$ has a unique best  symmetric nonnegative rank-one approximation.
\end{cor} 

Let $X \subset \CC^n$ be a complex variety. For $x \in X$ and $u \notin X$, let $d_u(x) = \sum_{i=1}^n (u(i) - x(i))^2$. The \textit{Euclidean distance degree} (ED degree) of $X$ is the number of  nonsingular critical points of $d_u$ for a generic $u$, and the \textit{ED discriminant} is the set of $u$ such that at least two critical points of $d_u$ coincide \cite{DrHOSTh15:fcm}. Hence Theorem~\ref{EDsym} shows that the ED discriminant of the cone over the Veronese variety (in both the real and complex case) is a hypersurface, and $D_{\operatorname{eig}}(T) = 0$ gives its defining equation.
\begin{example}
Let $T = [T_{ijk}] \in \mathsf{S}^3(\RR^2)$. Then $\psi_T(\lambda)$ is the resultant of the polynomials
\[
\begin{cases}
F_0 = T_{111}x^2 + 2T_{112}xy + T_{122}y^2 - \lambda xz, \\
F_1 = T_{112}x^2 + 2T_{122}xy + T_{222}y^2 - \lambda yz, \\
F_2 = x^2 + y^2 - z^2.
\end{cases}
\]
Consider the Jacobian determinant $J$ of $F_0, F_1, F_2$. Then\par
\vspace*{-2ex}
{\footnotesize
\begin{align*}
J &=  \det
\begin{bmatrix}
\partial F_0/\partial x & \partial F_0/\partial y & \partial F_0/\partial z \\
\partial F_1/\partial x & \partial F_1/\partial y & \partial F_1/\partial z \\
\partial F_2/\partial x & \partial F_2/\partial y & \partial F_2/\partial z
\end{bmatrix}\\
&= (8T^2_{112} - 8T_{111}T_{122} - 2 \lambda^2)x^2z + 4T_{112}\lambda y^3 \\
&\qquad + (8T_{122}^2 - 8T_{222}T_{112} - 2\lambda^2) y^2z + 4T_{122}\lambda x^3 \\
&\qquad + (4T_{111}\lambda - 8T_{122}\lambda) xy^2 + (4T_{122}\lambda + 4T_{111}\lambda) xz^2 \\
&\qquad + (4T_{112}\lambda + 4T_{222}\lambda) yz^2 + (4T_{222}\lambda - 8T_{112}\lambda) x^2y \\
&\qquad - 2\lambda^2z^3 + (8T_{112}T_{122} - 8T_{111}T_{222}) xyz,\\
\frac{\partial J}{\partial x} &= 12T_{122}\lambda x^2 + (4T_{111}\lambda - 8T_{122}\lambda) y^2 \\
& \qquad + (4T_{111}\lambda + 4T_{122}\lambda) z^2 + (8T_{222}\lambda - 16T_{112}\lambda) xy \\
& \qquad + (16T_{112}^2 - 4\lambda^2 - 16T_{111}T_{122}) xz \\
& \qquad + (8T_{112} T_{122} - 8T_{111} T_{222}) yz, \\
\frac{\partial J}{\partial y} &= (4T_{222}\lambda - 8T_{112}\lambda) x^2 + 12T_{112}\lambda y^2 \\
& \qquad + (4T_{112}\lambda + 4T_{222}\lambda) z^2 + (8 T_{111}\lambda - 16T_{122}\lambda) xy \\
& \qquad + (8T_{112}T_{122} - 8T_{111}T_{222}) xz \\
& \qquad + (16T_{122}^2 - 4\lambda^2 - 16T_{112}T_{222}) yz, \\
\frac{\partial J}{\partial z} &= (8T_{112}^2 - 8T_{111}T_{122} - 2\lambda^2) x^2 + (8T_{122}\lambda + 8T_{111}\lambda) xz \\
& \qquad + (8T_{122}^2- 8T_{222} T_{112} - 2\lambda^2) y^2 - 6\lambda^2 z^2 \\
& \qquad + (8T_{112}T_{122} - 8T_{111}T_{222}) xy + (8T_{112}\lambda + 8T_{222}\lambda) yz .
\end{align*}}%
By Salmon's formula \cite{CLO05}, $\psi_T(\lambda) = \frac{1}{512} \det (G)$, where $G$ is defined by \eqref{eq:G}.
Thus $\psi_T(\lambda) = p_2\lambda^6 + p_4\lambda^4 + p_6\lambda^2 + p_8$ for some homogeneous polynomials $p_m$ of degree $m$ in $T_{ijk}$, $m =2,4,6,8$. See also \cite{CS13:laa, LQZ13:cms}. Therefore $D_{\operatorname{eig}}(T)$ is the determinant of some $11 \times 11$ matrix in $T_{ijk}$.

For a generic $T\in \mathsf{S}^3(\RR^2)$, $\psi_T(\lambda) = c(\lambda^2 - \gamma_1)(\lambda^2 - \gamma_2)(\lambda^2 - \gamma_3)$ for some $c \in \CC$ and distinct $\gamma_i \in \CC$, and so $D_{\operatorname{eig}}(T) \neq 0$.

For $T \in H_{\operatorname{disc}}$, $\psi_T(\lambda)$ has multiple roots. For a specific example, let $S \in \mathsf{S}^3(\RR^2)$ be defined by $S_{111} = S_{222} = 1$ and set other  $S_{ijk} = 0$. Then $D_{\operatorname{eig}}(S) = 0$, implying that $S$ has at least one nonsimple eigenpairs. In fact, $\psi_{S}(\lambda) = (\lambda + 1)^2(\lambda - 1)^2(2\lambda^2 - 1)$ and so $S$ has two eigenvectors $(1, 0)$, $(0, 1)$ corresponding to the eigenvalue $1$, and two eigenvectors $(-1, 0)$, $(0, -1)$ corresponding to the eigenvalue $-1$. Note that $S$ is, up to a change of coordinates, the same example mentioned  at the beginning of this section, i.e., $S = x^{\otimes 3} + y^{\otimes 3}$ has two best rank-one approximations $x^{\otimes 3}$ and $y^{\otimes 3}$.
\end{example}

\section{Uniqueness of best rank-one approximations for real tensors}\label{sec:gen-unique}

In this section, $V$ and $W$, with or without subscripts, would generally denote real and complex vector spaces respectively.

Let $W_1, \dots, W_d$ be complex vector spaces. For $T \in W_1 \otimes \cdots \otimes W_d$, $u_i \in W_i$, and $\alpha_i \in \CC$, we denote by $\varphi_T(\lambda)$ the resultant of the following homogeneous polynomial equations
\begin{equation}\label{eq:rescha}
\begin{cases}
\alpha_i \langle T, u_1 \otimes \dots \otimes \widehat{u_i} \otimes \dots \otimes u_d \rangle = \lambda (\prod_{j \neq i} \alpha_j) u_i, \\
\langle u_i, u_i \rangle = \alpha_i^2,
\end{cases}
\end{equation}
for $i =1,\dots,d$.
Again by standard theory of resultants \cite{CLO05, GKZ94}, $\varphi_T(\lambda)$ vanishes if and only if \eqref{eq:rescha} has a nontrivial solution, and we obtain the following analogue of Definition~\ref{def:Echar}.
\begin{defi} $\varphi_T(\lambda)$ is called the \textit{singular characteristic polynomial} of $T \in W_1 \otimes \cdots \otimes W_d$.
\end{defi}
Clearly the roots of $\varphi_T(\lambda)$ are the normalized singular values of $T$. We also have an analogue of Definition~\ref{def:eigpair}.

\begin{defi}
Let $T \in W_1 \otimes \cdots \otimes W_d$. Two normalized singular pairs $(\lambda, u_1, \dots, u_d)$ and $(\mu, v_1, \dots, v_d)$ of $T$ are called \textit{equivalent} if $(\lambda, u_1, \dots, u_d) = (\mu, v_1, \dots, v_d)$, or $(-1)^{d-2}\lambda = \mu$ and $u_i = - v_i$ for $i = 1, \dots, d$. A normalized singular value $\lambda$ of $T$ is said to be \textit{simple} if it has only one corresponding normalized singular pair up to equivalence.
\end{defi}
For real vector spaces $V_1$, $\dots$, $V_d$, and $T \in V_1 \otimes \cdots \otimes V_d$, normalized singular pairs are invariant under the product of orthogonal groups $\O(n_1) \times \cdots \times \O(n_d)$.

It follows from \cite{FrieS16:banach} that the subset $X \subseteq V_1 \otimes \cdots \otimes V_d$ consisting of tensors  without unique best rank-one approximations is contained in a hypersurface. We will show that this can be strengthened to an \textit{algebraic} hypersurface.

\begin{proposition}\label{thm:3tenhy}
The following subset is  an algebraic hypersurface in  $ V_1 \otimes \cdots \otimes V_d$,
\begin{multline*}
X := \{ T \in V_1 \otimes \cdots \otimes V_d :  T \; \text{has non-unique}\\
\text{best rank-one approximations}\}. 
\end{multline*}
\end{proposition}

\begin{proof}
By Lemma~\ref{le:rel}, $X$ comprises tensors $T$ for which $\lVert T \rVert_\sigma$ is not  a simple normalized singular value. Let $d=3$ for notational simplicity. Let $T \in X$. Then there exist some $v_1, v_2, v_3$ with $\norm{v_i} = 1$ and $\{u_{1,1}, u_{2,1}, u_{3,1}\} \neq \{v_1, v_2, v_3\}$ with $\norm{u_{i,1}} = 1$ such that
\[
\langle T, u_{1,1} \otimes u_{2,1} \otimes u_{3,1}\rangle = \lVert T \rVert_\sigma = \langle T, v_1 \otimes v_2 \otimes v_3 \rangle.
\]
For each $i=1,2,3$, extend $u_{i,1}$ to an orthonormal basis $\{u_{i,1}, \dots, u_{i,n_i} \}$ of $V_i$. By an action of $\O(n_1) \times \O(n_2) \times \O(n_3)$ on $V_1 \otimes V_2 \otimes V_3$, we may assume that $v_i = \cos \theta_i u_{i,1} + \sin \theta_i u_{i,2}$. Let $T_{ijk} = \langle T, u_{1,i} \otimes u_{2,j} \otimes u_{3,k}\rangle$. Then we have\par
\vspace*{-2ex}
{\scriptsize
\begin{equation}\label{eq:tenpara}
\begin{cases}
T_{111} = \lVert T \rVert_\sigma, \\
T_{i11} = T_{1i1} = T_{11i} = 0, \\
T_{111}\cos \theta_2 \cos \theta_3 + T_{122}\sin \theta_2 \sin \theta_3 = T_{111} \cos \theta_1, \\
T_{212} \cos \theta_2 \sin \theta_3 + T_{221} \sin \theta_2 \cos \theta_3 + T_{222} \sin \theta_2 \sin \theta_3 = T_{111} \sin \theta_1, \\
T_{j12} \cos \theta_2 \sin \theta_3 + T_{j21} \sin \theta_2 \cos \theta_3
+ T_{j22} \sin \theta_2 \sin \theta_3 = 0, \\
T_{111}\cos \theta_1 \cos \theta_3 + T_{212}\sin \theta_1 \sin \theta_3 = T_{111} \cos \theta_2, \\
T_{122} \cos \theta_1 \sin \theta_3 + T_{221} \sin \theta_1 \cos \theta_3 + T_{222} \sin \theta_1 \sin \theta_3 = T_{111} \sin \theta_2, \\
T_{1j2} \cos \theta_1 \sin \theta_3 + T_{2j1} \sin \theta_1 \cos \theta_3
+ T_{2j2} \sin \theta_1 \sin \theta_3 = 0, \\
T_{111}\cos \theta_1 \cos \theta_2 + T_{221}\sin \theta_1 \sin \theta_2 = T_{111} \cos \theta_3, \\
T_{122} \cos \theta_1 \sin \theta_2 + T_{212} \sin \theta_1 \cos \theta_2 + T_{222} \sin \theta_1 \sin \theta_2 = T_{111} \sin \theta_3, \\
T_{12j} \cos \theta_1 \sin \theta_2 + T_{21j} \sin \theta_1 \cos \theta_2
+ T_{22j} \sin \theta_1 \sin \theta_2 = 0, 
\end{cases}
\end{equation}}%
for $i \neq 1$ and $j > 2$. By eliminating the parameter $\theta$, we obtain a system of polynomial equations that the $T_{ijk}$'s satisfy.

Let $J$ be the \textit{incidence variety} in $V_1\otimes V_2 \otimes V_3 \times \O(n_1) \times \O(n_2) \times \O(n_3)$, i.e., for each $(T, g_1, g_2, g_3) \in J$ where $g_i = [u_{i,1}, \dots, u_{i,n_i}] \in \O(n_i)$, there is some $(\theta_1, \theta_2, \theta_3)$ such that the $T_{ijk}$'s satisfy \eqref{eq:tenpara}. Define the projections
\begin{equation}
\xymatrix{
& J \ar[ld]^{\pi_{1}} \ar[rd]_{\pi_{2}} &\\
V_1 \otimes V_2 \otimes V_3 & & \O(n_1) \times \O(n_2) \times \O(n_3) }
\end{equation}
by $\pi_1 (T, g_1, g_2, g_3) = T$ and $\pi_2(T, g_1, g_2, g_3) = (g_1, g_2, g_3)$. Since $\lVert T \rVert_\sigma$ is a root of $\varphi_T(\lambda)$, $\lVert T \rVert_\sigma$ and its normalized singular vector tuples depend algebraically on $T$, implying that $J$ is an algebraic variety. $\lVert T \rVert_\sigma$ is simple if and only if $T$ is in the image of $\pi_1$, i.e., $X = \pi_1 (J)$.

Define $T' \in V_1\otimes V_2 \otimes V_3$ by $T'_{111} = T'_{222} = 1$ and set all other terms $T'_{ijk} = 0$. Then $T'$ has two normalized singular vector tuples corresponding to its normalized singular value$\lVert T \rVert_\sigma$. So $T' \in \pi_1(J)$. Since $T'$ has a finite number of singular pairs, a generic $T \in \pi_1(J)$ must also have a finite number of singular pairs by semicontinuity. So $\dim \pi_1^{-1}(T) = \dim \O(n_1-2) + \dim \O(n_2-2) + \dim \O(n_3-2)$ for a generic $T \in \pi_1(J)$, and $\dim X = \dim \pi_1(J) = \dim J - \dim \O(n_1-2) - \dim \O(n_2-2) - \dim \O(n_3-2)$.

Since $\pi_2$ is a dominant morphism, and the dimension of a generic fiber $\pi_2^{-1}(g_1, g_2, g_3)$ is $\dim V_1 \otimes V_2 \otimes V_3 - 2(n_1 + n_2 +n_3) +8$, it follows that $\dim J = \dim V_1 \otimes V_2 \otimes V_3 - 2(n_1 + n_2 +n_3) + 8 + \dim \O(n_1) + \dim \O(n_2) + \dim \O(n_3)$. Therefore the codimension of $X$ is $1$.
\end{proof}

We will show that normalized singular vector tuples of a generic tensor are distinct, a result that we will need later but is also of independent interest.
\begin{proposition}\label{prop:unisvt}
Let $W_1,\dots, W_d$ be vector spaces over $\CC$. A generic $T \in W_1\otimes \cdots \otimes W_d$  has distinct equivalence classes of normalized singular pairs.
\end{proposition}

Our proof of Proposition~\ref{prop:unisvt} will rely on the next three intermediate results. The first required result is a `Bertini-type' theorem introduced in \cite{FrieO14:fcm}.
\begin{theorem}[Friedland--Ottaviani]\label{thm:vb}
Let $E$ be a vector bundle on a smooth variety $B$. Let $S \subseteq H^0(B, E)$ generate $E$. If $\rank E > \dim B$, then the zero locus of a generic $\zeta \in S$ is empty.
\end{theorem}
\begin{lemma}\label{le:open}
Let $T \in W_1\otimes \cdots \otimes W_d$ be generic and let $(u_1, \dots, u_d)$ be a normalized singular vector tuple of $T$. If $v_d$ is not a scalar multiple of $u_d$, then $(u_1, \dots, u_{d-1},v_d)$ is not a normalized singular vector tuple of $T$.
\end{lemma}
\begin{proof}
Suppose $\lambda u_d = \langle T, u_1 \otimes \cdots \otimes u_{d-1} \rangle = \mu v_d$ for some $v_d$ not a scalar multiple of $u_d$. Then $\lambda$ or $\mu$ must be $0$, contradicting the fact that $0$ cannot be a singular value of a generic $T$ \cite[Theorem~1]{FrieO14:fcm}.
\end{proof}

\begin{lemma}\label{lem:gen}
Let $u_i, v_i, w_i \in W_i$ with $\langle u_i, u_i \rangle = \langle v_i, v_i \rangle = 1$, $i=1,2, 3$. For $x \in W_i$, we write $[x]_i$ for the corresponding element in the quotient space $W_i/\operatorname{span}( u_i)$. Suppose $u_i= v_i$ for at most one $i$. Then 
\begin{enumerate}[\upshape (i)]
\item the system of linear equations
\begin{equation}\label{eq:FOs}
\begin{cases}
\langle T, u_2 \otimes u_3 \rangle = \langle T, v_1 \otimes v_2 \otimes v_3 \rangle u_1 + w_1, \\
\langle T, u_1\otimes u_3 \rangle = w_2, \\
\langle T, u_1 \otimes u_2 \rangle = w_3,
\end{cases}
\end{equation}
has a solution $T \in W_1\otimes W_2 \otimes W_3$ if and only if $\langle u_2, w_2 \rangle = \langle u_3, w_3 \rangle$;
\item the system of linear equations
\begin{equation}\label{eq:FOq}
\begin{cases}
\langle T, u_2 \otimes u_3 \rangle = \langle T, v_1 \otimes v_2 \otimes v_3 \rangle u_1 + w_1, \\
[\langle T, u_1\otimes u_3 \rangle]_2 = [w_2]_2, \\
[\langle T, u_1 \otimes u_2 \rangle]_3 = [w_3]_3,
\end{cases}
\end{equation}
always has a solution $T \in W_1\otimes W_2 \otimes W_3$.
\end{enumerate}
\end{lemma}

\begin{proof} Note that the variables in these linear equations are $T_{ijk}$'s, the coordinates of $T$.
\begin{enumerate}[\upshape (i)]
\item Let $A$ be the coefficient matrix in \eqref{eq:FOs} and $b$ be the right-hand side. The system has a solution if and only if $A$ and the augmented matrix $[A \mid b]$ have the same rank, i.e., if there is some $x_i \in W_i$, $i =1,2,3$, such that $x_1 \otimes u_2 \otimes u_3 + u_1 \otimes x_2 \otimes u_3 + u_1 \otimes u_2 \otimes x_3 - \langle x_1, u_1 \rangle \cdot v_1 \otimes v_2 \otimes v_3 = 0$, then $\langle x_1, w_1 \rangle + \langle x_2, w_2 \rangle + \langle x_3, w_3 \rangle = 0$. Since $x_1 \otimes u_2 \otimes u_3 + u_1 \otimes x_2 \otimes u_3 + u_1 \otimes u_2 \otimes x_3 - \langle x_1, u_1 \rangle \cdot v_1 \otimes v_2 \otimes v_3 = 0 $ if and only if $x_1 = 0$, $x_2 = \alpha u_2$, $x_3 = -\alpha u_3$ or $x_1 = 0$, $x_2 = -\alpha u_2$, $x_3 = \alpha u_3$ for some $\alpha$, the system \eqref{eq:FOs} has a solution if and only if $\langle u_2, w_2 \rangle = \langle u_3, w_3 \rangle$.
\item The system \eqref{eq:FOq} has a solution if and only if $\langle u_2, w_2+t_2u_2 \rangle = \langle u_3, w_3+t_3u_3 \rangle$ for some $t_2, t_3 \in \CC$. Choose any $t_2, t_3$ such that $t_3-t_2 = \langle u_2, w_2 \rangle - \langle u_3, w_3 \rangle$. \qedhere
\end{enumerate}
\end{proof}

\begin{proof}[Proof of Proposition~\ref{prop:unisvt}]
Let $d = 3$ for notational convenience. For $i =1,2,3$, let $C_i = \{u_i \in W_i: \langle u_i, u_i \rangle = 1 \}$, $F_i$ be the trivial vector bundle on $C_i$ with fiber isomorphic to $W_i$, $L_i$ be the tautological line bundle on $C_i$, and $Q_i$ be the quotient bundle $F_i/L_i$ on $C_i$. Consider the exact sequence of vector bundles
\[
0 \to L_i \to F_i \to Q_i \to 0
\]
over $C_i$. Let $M = C_1 \times C_2 \times C_3$. We will need to discuss vector bundles over $M \times M$ and for clarity, we distinguish the two copies of $M$. So we write  $M_1 \times M_2$ where  $M_1 = M_2 = M$. Let $\pi_{i, j} : M_i \rightarrow C_j$ be the natural projection for $i = 1, 2$ and $j = 1, 2, 3$. Let $p_i: M_1 \times M_2 \rightarrow M_i$ be the natural projection for $i = 1, 2$. Then we have the following diagram:\par
\vspace*{-1ex}
{\footnotesize
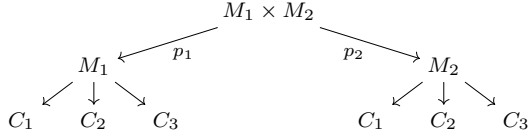
\[
\begin{tikzcd}[column sep=small,row sep=small]
& & & M_1 \times M_2 \arrow{dll}{p_{1}} \arrow{drr}[swap]{p_{2}} & & &\\
& M_1 \arrow{dl} \arrow{d} \arrow{dr} & & & & M_2 \arrow{dl} \arrow{d} \arrow{dr} &\\
C_1 & C_2 & C_3 &  & C_1 & C_2 & C_3
\end{tikzcd}.
\]}%
Consider the vector bundle on $M_1 \times M_2$,
\[
\widetilde{E} = \biggl(\bigoplus\limits_{j=1}^3 p_1^* \pi_{1,j}^* (Q_j)\biggr) \oplus p_2^*\pi_{2,1}^*(F_1) \oplus \biggl(\bigoplus \limits_{j=2}^3 p_2^*\pi_{2,j}^*(Q_j)\biggr),
\]
where $f^*$ denotes the pullback induced by $f$. Let
\begin{multline*}
X_i = \{(v_1,v_2,v_3,u_1,u_2,u_3)\in M_1 \times M_2:\\
 u_j = v_j \;\text{for all}\;  j\neq i \}.
\end{multline*}
By Lemma~\ref{le:open},  to study the behavior of normalized singular pairs of a generic tensor, we need only consider the following open subset of the affine variety $M_1 \times M_2$,
\[
B = M_1 \times M_2 \setminus (X_1 \cup X_2 \cup X_3),
\]
and its corresponding vector bundle $E = \widetilde{E} \vert_B$ over the base space $B$. Then
\[
\rank {E} = 2 \sum_{i=1}^3 \dim W_i - 5 > \dim B = 2 \sum_{i=1}^3 \dim W_i - 6.
\]
So the inequality in  Theorem~\ref{thm:vb} holds for our choice of $E$ and $B$.
Now let
\begin{align*}
S &= \{s \in H^0(B, E): s(v_1, v_2, v_3, u_1, u_2, u_3) \\
& =([\langle T, v_2\otimes v_3\rangle]_1, [\langle T, v_1\otimes v_3\rangle]_2, [\langle T, v_1\otimes v_2\rangle]_3, \\
& \qquad\qquad  \langle T, u_2\otimes u_3\rangle - \langle T, v_1 \otimes v_2 \otimes v_3 \rangle u_1, \\
& \qquad\qquad\qquad [\langle T, u_1\otimes u_3 \rangle]_2, [\langle T, u_1 \otimes u_2 \rangle]_3) \}.
\end{align*}
By Lemma~\ref{lem:gen} and \cite[Lemma~8]{FrieO14:fcm}, $S$ generates $E$. By Theorem~\ref{thm:vb}, a generic section of $E$ does not vanish on $B$, implying tha each normalized singular value of a generic tensor is distinct and simple.
\end{proof}

Let  $D_{\operatorname{sing}}(T)$ be the \textit{singular discriminant}, the resultant of the singular characteristic polynomial $\varphi_T$ and its derivative $\varphi'_T$. Since a generic $T$ has distinct equivalence classes of normalized singular pairs, $\varphi_T$ has simple roots, and so  $D_{\operatorname{sing}}(T)$ does not vanish identically. As $D_{\operatorname{sing}}(T)$ vanishes on $X$, the hypersurface defined in Theorem~\ref{thm:3tenhy},  $D_{\operatorname{sing}}(T) = 0$ indeed defines a hypersurface in $W_1 \otimes \cdots \otimes W_d$. 
Note that $X$ is a union of some components of the real points of $X_{\operatorname{disc}}$. Finally, we arrive at our main result of this section, singular value analogues of Theorem~\ref{EDsym} and Corollaries~\ref{cor:EDsym} and \ref{cor:EDnonneg}.
\begin{theorem}\label{main}
$D_{\operatorname{sing}}(T) = 0$ is a defining equation of the hypersurface
\begin{multline*}
X_{\operatorname{disc}} := \{ T \in W_1 \otimes \dots \otimes W_d :  T \; \text{has a}\\
\text{non-simple normalized singular value}\}. 
\end{multline*}
\end{theorem}
In the following, let $V_i$ be a real vector space and $W_i = V_i \otimes_\mathbb{R} \mathbb{C}$ be its complexification, $i =1,\dots,d$.
\begin{cor}\label{cor:main}
Let $ T \in V_1 \otimes \dots \otimes V_d$ be real. If $D_{\operatorname{sing}}(T) \neq 0$, then $T$ has a unique best rank-one approximation.
\end{cor}
We deduce the following analogue for nonnegative tensors from Theorem~\ref{thm:nonneg} and Corollary~\ref{cor:main}.
\begin{cor}
Let $T \in V_1 \otimes \dots \otimes V_d$ be nonnegative. If $D_{\operatorname{sing}}(T) \neq 0$, then $T$ has a unique best nonnegative rank-one approximation.
\end{cor} 

Theorem~\ref{main} shows that the ED discriminant $X_{\operatorname{disc}}$ of the cone over the Segre variety $\mathbb{P}W_1 \times \cdots \times \mathbb{P}W_d$ is a hypersurface when $d \ge 3$, and $D_{\operatorname{sing}}(T) = 0$ gives its defining equation. The discussion before Theorem~\ref{main} shows that the set of real points of $X_{\operatorname{disc}}$ is a real hypersurface. It is interesting to note that when $d = 2$, i.e., the matrix case, the set of real points of the ED discriminant of the Segre variety $\mathbb{P}W_1 \times \mathbb{P}W_2$ has codimension $2$ \cite[Example~7.6]{DrHOSTh15:fcm}.


\appendix

We use semirings and semimodules instead of rings and modules to construct tensor products of cones in order to give nonnegative tensors an algebraic description and state our results in a more general setting. A semimodule over a semiring is essentially the same notion as a vector space over a field, except that the field of scalars is now replaced by a semiring of scalars like the nonnegative reals. The nonnegative reals do not form a field or even a ring since they do not have additive inverses, but aside from this, $\mathbb{R}_+$ has all the properties of scalars that makes the notion of a vector space so useful in engineering.

\begin{defi}
A semiring $R$ is a set equipped with binary operations $+$ and $\cdot$ such that
\begin{itemize}
\item $(R, +)$ is a commutative monoid with identity element $0$;
\item $(R, \cdot)$ is a monoid with identity element 1;
\item Multiplication left and right distributes over addition:
\begin{align*}
a\cdot (b + c) = (a\cdot b) + (a\cdot c),\\
(a + b)\cdot c = (a\cdot c) + (b\cdot c);
\end{align*}
\item Multiplication by 0 annihilates $R$:
\[
0\cdot a = a\cdot 0 = 0.
\]
\end{itemize}
\end{defi}

\begin{defi}
A commutative semiring is a semiring whose multiplication is commutative.
\end{defi}

\begin{defi}
A semimodule $M$ over a commutative semiring $R$ is a commutative monoid $(M, +)$ and an operation $\cdot: R \times M \to M$ such that for all $r, s$ in $R$ and $x, y \in M$, we have:
\begin{align*}
r \cdot ( x + y ) &= r \cdot x + r \cdot y, \\
( r + s ) \cdot x &= r \cdot x + s \cdot x, \\
( rs ) \cdot x &= r \cdot ( s \cdot x ), \\
1_R \cdot x &= x.
\end{align*}
\end{defi}
In our context, the set of nonnegative real numbers $\RR_+$ is a commutative semiring and the set of nonnegative tensors is a semimodule over $\RR_+$.

\section*{Acknowledgment}

YQ thanks Giorgio Ottaviani for careful reading and very helpful advice, especially for suggesting that we use an argument in \cite{FrieO14:fcm} for Proposition~\ref{prop:unisvt}. YQ also thanks Emil Horobe{\c t} for helpful discussion. LHL thanks Shmuel Friedland for a pointer to Rademacher Theorem.  We thank Yuning Yang for pointing out an error in an earlier version. Special thanks to the three anonymous reviewers for their useful suggestions.

\bibliographystyle{ieeetranS}

\end{document}